\documentclass[12pt]{elsarticle}
\usepackage{graphicx}
\usepackage{amssymb}    
\usepackage{epstopdf}
\usepackage{comment}
\usepackage{color}
\usepackage{hyperref}
\usepackage{algorithm,algorithmic}
\usepackage[final]{pdfpages}
\usepackage{frcursive}
\usepackage{bbm}
\usepackage{mathtools}
\usepackage{tabularx}
\usepackage{algorithm,algorithmic}

\usepackage{yfonts}
\usepackage{lscape}
\usepackage{epsfig}
\usepackage{fullpage}
\usepackage[titletoc]{appendix}
\usepackage{fancyhdr}
\usepackage{amsmath,amssymb,xspace}

\usepackage{tabularx}
\usepackage{mathrsfs}

\newfont{\bbb}{msbm10 scaled\magstep1}

\newtheorem{rem}{Remark}[section]

\let \leq \leqslant
\let \geq \geqslant

\let \epsilon \varepsilon

{ \par \medskip \par
  \noindent \textit{\textbf{Demonstration\/}} : }{\null \hfill $\Box$ \par }
 
\newcommand{\R} {\ensuremath{\mathbb{R}}}
{

\newcommand{\N} {\ensuremath{\mathbb{N}}}

\makeatletter
\newcommand{\doublewidetilde}[1]{{%
  \mathpalette\double@widetilde{#1}%
}}
\newcommand{\double@widetilde}[2]{%
  \sbox\z@{$\m@th#1\widetilde{#2}$}%
  \ht\z@=.9\ht\z@
  \widetilde{\box\z@}%
}

\begin{document}

\begin{frontmatter}


\title{Schwarz Waveform Relaxation Physics-Informed Neural Networks for Solving Advection-Diffusion-Reaction Equations}

\author[carl,crm]{Emmanuel LORIN}
\ead{elorin@math.carleton.ca}
\author[ucsb]{Xu YANG}
\ead{xuyang@math.ucsb.edu}

\address[carl]{School of Mathematics and Statistics, Carleton University, Ottawa, Canada, K1S 5B6}
\address[crm]{Centre de Recherches Math\'{e}matiques, Universit\'{e} de Montr\'{e}al, Montr\'{e}al, Canada, H3T~1J4}
\address[ucsb]{Department of Mathematics, University of California, Santa Barbara, CA 93106, USA}

\begin{abstract}
This paper develops a physics-informed neural network (PINN) based on the Schwarz waveform relaxation (SWR) method for solving local and nonlocal advection-diffusion-reaction equations. Specifically, we derive the formulation by constructing subdomain-dependent local solutions by minimizing local loss functions, allowing the decomposition of the training process in different domains in an embarrassingly parallel procedure. Provided the convergence of PINN, the overall proposed algorithm is convergent. By constructing local solutions, one can, in particular, adapt the depth of the deep neural networks, depending on the solution's spectral space and time complexity in each subdomain. We present some numerical experiments based on classical and Robin-SWR to illustrate the performance and comment on the convergence of the proposed method.

\end{abstract}



\end{frontmatter}


\section{Introduction}


This paper focuses on the derivation and analysis of a Neural-Network (NN) based Schwarz Waveform Relaxation (SWR) Domain Decomposition Method (DDM) for solving partial differential equations (PDE) in parallel. We will focus in this paper on a simple diffusion-advection-reaction equation. Still, the proposed strategy applies to any other evolution (in particular wave-like) equations, for which the convergence of SWR is proven \cite{XA1,XA2,XA3,XA4,halpern2,halpern3,GanderHalpernNataf,DDMN,GanderHalpernNataf2,GanderReview,DoleanBook}. We derive a combined DDM-SWR and Physics-Informed Neural Network (PINN) method for solving local and nonlocal diffusion-advection-reaction equations. The latter was developed by Karniadakis {\it et al.} \cite{pinns,pinns2,pinns3} and is a general strategy in scientific machine learning for solving PDE using deep neural networks via the minimization of well-designed loss functions.  Notice that in \cite{pinns4}, was also proposed a more direct DDM for solving PDE. Interestingly, both (\cite{pinns4} and the one presented here) methods could actually be combined; this was however not tested in this paper.  Let us also mention a recent paper \cite{Heinlein2021}, where a combination of Schwarz DDM with NN-based solvers is proposed for stationary PDE. Beyond the derivation of the SWR-NN method, this paper's objective is to exhibit some fundamental properties that make this methodology very promising. The general principle is to solve Initial Boundary Value Problems (IBVPs) by constructing local solutions (subdomain-dependent) obtained by minimizing local loss functions. The overall strategy is convergent (provided that the PINN method is convergent)  and allows, in particular, to locally decompose the training process in different subdomains within an embarrassingly parallel procedure. The construction of local solutions also allows to locally adapt the depth of the deep neural network, depending on the solution's spectral space and time complexity in each subdomain.

In this paper, we will primarily focus on the derivation aspects and will not necessarily detail all the computational aspects, particularly regarding the selection of the training points. This will, however, be specified in the Numerical Experiment section. For convenience, we shall recall some basic aspects about PINNs, neural networks and SWR method for evolution equations, which shall be used in the paper later.

\subsection{Basics on PINNs}
Let us recall the principle of PINNs for solving, {\it e.g.}, an evolution PDE over $\Omega\times [0,T]$,
\begin{eqnarray}\label{pde0}
\left\{
\begin{array}{lcl}
\partial_t u + Pu & = & f \, \, \hbox{ in } \Omega\times[0,T] \\
Mu & = & 0 \, \, \hbox{ in } \Gamma\times[0,T]  \\
u(\cdot,0) & = & u_0  \, \, \hbox{ in } \Omega
\end{array}
\right. \, ,
\end{eqnarray}
where i) $P$ is a differential operator in space, and $M$ a differential or algebraic boundary operator over the domain boundary $\Gamma$; ii) are $f$ and $u_0$ are imposed functions. The PINN approach, which generalizes the DE-solver from Lagaris \cite{lagaris} consists in parameterizing (by say ${\boldsymbol w}$) a NN, $N({\boldsymbol w},{\boldsymbol x},t)$ approximating the solution to \eqref{pde0}, by minimizing (a discrete version of) the  following loss function
\begin{eqnarray*}
\mathcal{L}({\boldsymbol w}) & = &\|(\partial_t + P)N({\boldsymbol w},\cdot,\cdot) -f\|_{L^2(\Omega\times (0,T))} + \lambda_1\|MN({\boldsymbol w},\cdot,\cdot)\|_{L^2(\Gamma\times (0,T))}\\
& &  + \lambda_2\|N({\boldsymbol w},\cdot,0) - u_0\|_{L^2(\Omega)} \, ,
\end{eqnarray*}
where $\lambda_{1,2}$ are some free positive parameters and ${\boldsymbol w} \in W \in \R^P$ for some large $P$, and where $\|\cdot\|_{L^2(\Omega\times (0,T))}$ (resp. $\|\cdot\|_{L^2(\Gamma\times (0,T))}$) denotes the $L^2-$norm over $\Omega\times(0,T)$ (resp. $\Gamma\times(0,T)$). Practically, the loss functions are constructed by estimating the values at a very large number of space$\&$time-training points $\{({\boldsymbol x}_j,t_n)\}_{j,n}$. Hence the $L^2$-norm are not exactly computed, but only approximated. Karniadakis and collaborators have developed numerous techniques to improve the efficiency and rate of convergence of PINN-algorithms for different types of PDE. We refer for instance to \cite{pinns,pinns2,pinns3} for details. Ultimately, the PINN strategy is to provide more efficient solvers than standard methods (finite difference, finite-volume, finite-elements, spectral methods, pseudospectral methods, etc) for high dimensional (stochastic or deterministic) PDEs. As far as we know, this is not clearly yet established, which justifies the active research in this field and the developed of new methods.

\subsection{Basics of Neural Networks}
We here recall the basics of neural networks. We denote the neural network $N({\boldsymbol w},{\boldsymbol x})$ with ${\boldsymbol x}=(x_1,\cdots,x_d) \in \Omega \subseteq \R^d$ and where we denote by ${\boldsymbol w}$ the unknown parameters. Neural networks usually read (for 1 hidden layer, machine learning)
\begin{eqnarray}\label{NN0}
\left.
\begin{array}{lcl}
N({\boldsymbol w},{\boldsymbol x}) & = & \sum_{i=1}^Hv_{i}\sigma_{i}\big(\sum_{j=1}^dw_{ij}x_j+u_i\big) \, ,
\end{array}
\right.
\end{eqnarray}
where $\{\sigma_i\}_{1\leq i \leq H}$ are the sigmoid transfer functions, and $H$ is the number of sigmoid units, $\{w_{ij}\}_{ij}$ are the weights and $\{u_i\}_i$ the bias. When considering several hidden layers (deep learning), we have then to compose functions of the form \eqref{NN0}, \cite{despres}. That is
\begin{eqnarray*}
N & = & \mathcal{N}_p \circ  \mathcal{N}_{p-1}\circ \cdots \mathcal{N}_1\circ  \mathcal{N}_0 \, ,
\end{eqnarray*}
where for $0\leq r\leq p$, $\mathcal{N}_r$ is defined from $\R^{a_r}$ (with $a_{r} \in \N$) to $\R^{a_{r+1}}$ by $\sigma_r(W_rX_r+b_r)$, $\sigma_r$ is an activation function, $X_r\in \R^{a_r}$ and where $(a_0,\cdots,a_{p+1})$ where $p+1$ layers are considered. The layer $r=0$ is the input layer and $r=p+1$ is the output layer, such that $a_0=a_{p+1}=m$. In fine, $N$ from $X\in \R^{m}$ to $\R^{m}$.

\subsection{Basics on SWR methods for evolution equations}
In this subsection, we recall the principle of SWR-DDM for solving evolution PDE. Consider a $d$-dimensional first order in time evolution partial differential equation $\partial_tu + Pu=f$ in the spatial domain $\Omega \subseteq \R^d$, and time domain $(0,T)$, where $P$ is a linear differential operator in space. The  initial data is denoted by $u_0$, and we impose, say, null Dirichlet boundary conditions on $\Gamma_{\textrm{ext}}=\partial\Omega$. We present the method for 2 subdomains, although in practice an arbitrary number of subdomains can be employed. We first split $\Omega$ into two open subdomains $\Omega_{\varepsilon}^{\pm}$, with or without overlap ($\Omega_{\varepsilon}^+\cap\Omega_{\varepsilon}^-=\emptyset$ or $\Omega_{\varepsilon}^+\cap\Omega_{\varepsilon}^-\neq \emptyset$), and $\epsilon \geq 0$. The SWR algorithm consists in iteratively solving IBVPs in $\Omega_{\varepsilon}^{\pm}\times (0,T)$, using transmission conditions at the subdomain interfaces $\Gamma_{\varepsilon}^{\pm}:=\partial\Omega_{\varepsilon}^{\pm}\cap\Omega_{\varepsilon}^{\mp}$. The imposed transmission conditions are established using the preceding Schwarz iteration data in the adjacent subdomain. That is, for $k\geq1$, and denoting $u^{\pm}$ the solution in $\Omega_{\varepsilon}^{\pm}$, we consider
\begin{eqnarray}\label{owrs}
\left\{
\begin{array}{lcl}
(\partial_t + P) u^{\pm, (k)} & = & f, \, \hbox{ in } \Omega_{\varepsilon}^{\pm} \times (0,T),\\
 u^{\pm, (k)}(\cdot,0) & =& u_0^{\pm},\, \hbox{ in } \Omega_{\varepsilon}^{\pm}, \\
\mathcal{T}_{\pm} u^{\pm, (k)} & = & \mathcal{T}_{\pm} u^{\mp, (k-1)}, \, \hbox{ on } \Gamma_{\varepsilon}^{\pm} \times (0,T),\\
u^{\pm, (k)} & = & 0, \, \hbox{ on } \Lambda_{\varepsilon}^{\pm} \times (0,T) \, ,
\end{array}
\right.
\end{eqnarray}
with a given initial guess $u^{\pm, (0)}$, where $\mathcal{T}^{\pm}$ denotes a boundary$/$transmission operator and where $\Lambda_{\varepsilon}^{\pm}:=\partial \Omega_{\varepsilon}^{\pm}\backslash \Gamma_{\varepsilon}^{\pm}$ are internal boundaries. Classical SWR (CSWR) method consists in taking $\mathcal{T}^{\pm}$ as the identity operator while Optimized-SWR (OSWR) method consists in taking $\mathcal{T}^{\pm} = \nabla_{{\boldsymbol n}^{\pm}} \pm \lambda^{\pm}_{\Gamma_{\varepsilon}}$ for some well chosen (optimized from the convergence rate point of view) $\lambda^{\pm}_{\Gamma^{\pm}_{\varepsilon}} \in \R^*$, and outward normal vector ${\boldsymbol n}^{\pm}$ to $\Gamma^{\pm}_{\varepsilon}$. The OSWR method is then a special case of Robin-SWR methods. In addition, in order to provide a faster convergence than CSWR, the OSWR method is often convergent even for non-overlapping DDM. The latter is hence of crucial interest from the computational complexity point. We refer to \cite{halpern3,GanderHalpernNataf2,GanderReview,DoleanBook,GanderHalpernNataf} for details.  The convergence criterion for the Schwarz DDM is typically given for any $0<t\leq T$, by
\begin{eqnarray}\label{CVTOT}
\big\| \hspace{0.2cm}\|u^{+,(k)}(\cdot,t)-u^{-,(k)}(\cdot,t)\|_{\infty;\overline{\Omega}_{\varepsilon}^+\cap\overline{\Omega}_{\varepsilon}^-}\big\|_{L^{2}(0,t)}\leq  \delta^{\textrm{Sc}},
\end{eqnarray}
 with $\delta^{\textrm{Sc}}$ small enough. When the convergence of the full iterative algorithm is obtained at Schwarz iteration $k^{\textrm{cvg}}$,  one gets the converged global solution 
$u^{\textrm{cvg}}:=u^{(k^{\textrm{cvg}})}$ in $\Omega$. The reconstructed solution $u$, is finally defined as $u_{|\Omega_{\varepsilon}^{\pm}}=u^{\pm,(k^{\textrm{cvg}})}$.
\subsection{Advection-diffusion-reaction equation}
Rather than considering a general situation, for which the rapid convergence of the SWR method and efficiency are not necessarily proven, we propose to focus on the advection-diffusion-reaction equation, for which both properties are established in \cite{halpern3} (see also \cite{halpern2,XA1,XA2,XA3,XA4} for the Schr\"odinger equation). Let us consider the following initial boundary-value problem: find the real function $u({\boldsymbol x},t)$ solution to the advection-diffusion-reaction equation on $ \R^d$, $d \geq 1$,
\begin{eqnarray}
\label{e1}
\left\{
\begin{array}{l}
\partial_t u = \nu({\boldsymbol x})\triangle u + {\boldsymbol a}({\boldsymbol x}) \cdot \nabla u + r({\boldsymbol x})u, \, {\boldsymbol x} \in  \R^d, \, t>0,\\
u({\boldsymbol x},0) = u_0({\boldsymbol x}), \, {\boldsymbol x} \in  \R^d,
\end{array}
\right.
\end{eqnarray} 
with initial condition $u_{0}$, and the real-valued space-dependent smooth reaction term $r$, advection vector ${\boldsymbol a}$ and diffusion $\nu$. \\
We recall from \cite{halpern3}, that for considering $\Omega=\R$ with constant coefficients in \eqref{e1}, and for $u_0 \in L^2(\Omega)$, $f\in L^2(0,T;L^2(\Omega))$, there exists a unique weak solution in $C(0,T;L^2(\Omega))\cap L^2(0,T;H^1(\Omega))$. Moreover, if $u_0\in H^2(\Omega)$ and $f\in L^2(0,T;L^2(\Omega))$, there exists a unique weak solution in $L^2(0,T;H^3(\Omega))\cap H^{3/2}(0,T;L^2(\Omega))$.

\subsection{Organization of the paper}
The rest of the paper is organized as follows. In Section \ref{sec:SWR-PINN}, we derive the combined SWR-PINN method, and some properties are proposed in Sections \ref{subsec:interest} and \ref{subsec:direct}. Section \ref{sec:numerics} is devoted to some numerical experiments illustrating the convergence of the overall SWR-PINN method. We make conclusive remarks in Section \ref{sec:conclusion}.

\section{SWR-PINN method}\label{sec:SWR-PINN}
In this section, we propose to combine PINN-based solvers with SWR-DDM to solve the advection-diffusion-reaction equation on a bounded domain $\Omega \subset \R^d$, imposing null Dirichlet boundary conditions at $\Gamma$. For the sake of simplicity of the presentation, the derivation is proposed for two subdomains; the extension to an arbitrary number of subdomains is straightforward.

\subsection{Derivation of the SWR-PINN method}\label{subsec:SWR-PINN}
The standard SWR method for two subdomains consists in solving the IVBP using the following algorithm
\begin{eqnarray}\label{oswr2}
\left\{
\begin{array}{lcl}
\partial_t u^{\pm, (k)} & = & \nu({\boldsymbol x})\triangle u^{\pm, (k)} +  {\boldsymbol a}({\boldsymbol x}) \cdot \nabla u^{\pm, (k)}+  r({\boldsymbol x})u^{\pm, (k)} , \, \hbox{ in } \Omega_{\varepsilon}^{\pm} \times (0,T),\\
 u^{\pm, (k)}(\cdot,0) & =& u_0^{\pm},\, \hbox{ in } \Omega_{\varepsilon}^{\pm}, \\
\mathcal{T}_{\pm} u^{\pm, (k)} & = & \mathcal{T}_{\pm} u^{\mp, (k-1)}, \, \hbox{ on } \Gamma_{\varepsilon}^{\pm} \times (0,T),\\
u^{\pm, (k)} & = & 0, \, \hbox{ on } \Lambda_{\varepsilon}^{\pm} \times (0,T) \, ,
\end{array}
\right.
\end{eqnarray}
where $\mathcal{T}_{\pm}$ is a boundary operator, and where we recall that $\Lambda_{\varepsilon}^{\pm}=\partial \Omega_{\varepsilon}^{\pm}\backslash \Gamma_{\varepsilon}^{\pm}$. The well-posedness and the convergence of this method and its rate of convergence were established in \cite{halpern3} for different types of transmission conditions. SWR algorithms can actually be reformulated as a fixed point methods (FPM), and their rate of convergence is hence determined by a contraction factor of the FPM. More specifically, it is proven in \cite{halpern2} that for
\begin{eqnarray*}
\partial_t u +a\partial_xu -\nu \partial_{xx}u +ru  = f \, \, \,  (x,t)\in \R^2\times  (0,T),
\end{eqnarray*}
where $f\in L^2(\Omega)$, $\nu>0$, $a,b$ in $\R$, the CSWR method is convergent and has a convergence rate $C_{\textrm{CSWR}}$ (contract factor), at least, given by 
\begin{eqnarray*}
C_{\textrm{CSWR}} &=& \exp\Big(-\cfrac{\varepsilon}{\nu}(\sqrt{a^2+4\nu r})\Big) \, .
\end{eqnarray*}
In fact, this can be refined to superlinear convergence rate $2/\sqrt{\pi}\int^{\infty}_{\varepsilon/\sqrt{\nu T}}e^{-s^2}ds$. For Robin-SWR methods, with transmission conditions $\partial_x \pm \lambda$,  the convergence rate is actually improved
\begin{eqnarray*}
C_{\textrm{Robin}} & = & \sup_{\omega \in \R}\Big|\cfrac{\sqrt{a^2+4\nu(r+{\tt i}\omega)}-\lambda}{\sqrt{a^2+4\nu(r+{\tt i}\omega)}+\lambda}\Big|\exp\Big(-\cfrac{\varepsilon}{\nu}(\sqrt{a^2+4\nu r})\Big) \, ,
\end{eqnarray*}
as $|\sqrt{a^2+4\nu(r+{\tt i}\omega)}-\lambda|/|\sqrt{a^2+4\nu(r+{\tt i}\omega)}+\lambda|<1$.  We notice in particular, that a crucial element for the rate of convergence of SWR methods, is the size of the overlapping zone. However, overlapping is not required for Robin-SWR methods to converge. \\
Rather than a standard approximation of \eqref{oswr2} using a finite elements$/$difference or pseudospectral methods \cite{XA4,halpern2,XA3}, we then propose to solve this system using a PINN method. We denote by $N({\boldsymbol w},{\boldsymbol x},t)$ the generic NN to optimize, where ${\boldsymbol w}$ denotes the unknown parameters. The SWR-NN hence consists in searching for an approximate solution to the SWR method by applying local PINN algorithms. That is, we now consider

\begin{eqnarray}\label{oswr3}
\left\{
\begin{array}{lcl}
\partial_t N^{\pm, (k)} & = & \nu({\boldsymbol x})\triangle N^{\pm, (k)} +  {\boldsymbol a}({\boldsymbol x}) \cdot \nabla N^{\pm, (k)}+ r({\boldsymbol x})N^{\pm, (k)}, \, \hbox{ in } \Omega_{\varepsilon}^{\pm} \times (0,T),\\
 N^{\pm, (k)}(\cdot,0) & =& u_0^{\pm},\, \hbox{ in } \Omega_{\varepsilon}^{\pm}, \\
\mathcal{T}_{\pm} N^{\pm, (k)} & = & \mathcal{T}_{\pm} N^{\mp, (k-1)}, \, \hbox{ on } \Gamma_{\varepsilon}^{\pm} \times (0,T),\\
N^{\pm, (k)} & = & 0, \, \hbox{ on } \Lambda_{\varepsilon}^{\pm} \times (0,T) \,.
\end{array}
\right.
\end{eqnarray}

\begin{rem}
For the CSWR method ($\mathcal{T}_{\pm}$ is the identity operator), it is proven in \cite{halpern3}, among many other well-posed results, that for $f\in L^2(0,T;H^1(\Omega_{\varepsilon}^{\pm}))$ and $u_0\in H^1(\Omega)$ with $u^{\pm,(0)}\in H^{3/4}(0,T)$ and some compatibility conditions, the algorithms \eqref{oswr2} is well-posed in $L^2(0,T;H^2(\Omega)\cap H^{1}(0,T;L^2(\Omega))$.\\
Let us now denote $e^{\pm,(k)}:=u^{\pm,(k)}-u$ and $\widetilde{e}^{\pm,(k)}:=N^{\pm,(k)}-u$. In Theorem 3.3 from \cite{halpern3} it is stated that for any $k\geq 0$ and for some positive constant $C>0$
\begin{eqnarray*}
\|(e^{+,(2k+1)},e^{-,(2k+1)})\|_{\mathcal{H}_{\varepsilon}} & \leq & C\exp\Big(-k\cfrac{\varepsilon(\sqrt{a^2+4\nu r})}{\nu}\Big)\|(u(\varepsilon,\cdot),u(-\varepsilon/2,\cdot))\|_{(_0H^{3/4}(0,T))^2} \, ,
\end{eqnarray*}
where we have denoted
\begin{eqnarray*}
\mathcal{H}_{\varepsilon} & := & L^2(0,T;H^2(\Omega_{\varepsilon}^{+})\cap H^{1}(0,T;L^2(\Omega_{\varepsilon}^{+})) \times L^2(0,T;H^2(\Omega_{\varepsilon}^{+})\cap H^{1}(0,T;L^2(\Omega_{\varepsilon}^{-})) \, .
\end{eqnarray*}
Now if we assume (convergence of the PINN-method) that the NN solution to \eqref{oswr3} is such that there exist $\eta(\overline{\boldsymbol w}^{\pm;(k)};\varepsilon)$
\begin{eqnarray*}
\|N^{\pm, (k)}-u^{\pm, (k)}\|_{\mathcal{H}_{\varepsilon}} & \leq  & \eta(\overline{\boldsymbol w}^{\pm;(k)};\varepsilon) \, ,
\end{eqnarray*}
then
\begin{eqnarray*}
\|\widetilde{e}^{\pm,(2k+1)}\|_{\mathcal{H}_{\varepsilon}} & \leq & \|e^{\pm,(2k+1)})\|_{\mathcal{H}_{\varepsilon}} + \|N^{\pm, (2k+1)}-u^{\pm, (2k+1)}\|_{\mathcal{H}_{\varepsilon}} \, .
\end{eqnarray*}
We then trivially deduce the convergence of the overall PINN-CSWR method. Similar conclusions can be reached for the PINN-OSWR algorithms. In particular in \cite{halpern3}, the OSWR method is shown to be convergent in $L^2(0,T;H^3(\Omega))\cap H^{3/2}(0,T;L^2(\Omega))$.
\end{rem}
 First, we use the standard technique to include the initial condition \cite{lagaris}, by searching for a trial network in the form
\begin{eqnarray*}
T({\boldsymbol w},{\boldsymbol x},t) & = & u_0({\boldsymbol x}) + tN({\boldsymbol w},{\boldsymbol x},t) \, .
\end{eqnarray*}
For the sake of simplicity of the notations, we will yet denote by $N$ the neural networks, which includes the contribution of the initial data.  Notice that this step is not essential, but allows to simplify the loss functions. Hence at each Schwarz iteration we solve \eqref{oswr3}, by minimization the following {\it two} local ``independent'' loss functions $\mathcal{L}^{\pm}$, for some positive parameters $\lambda_{\textrm{Int}}^{\pm}$, $\lambda_{\textrm{Ext}}^{\pm}$ and where we have denoted ${\boldsymbol w}=({\boldsymbol w}^-,{\boldsymbol w}^+)$. In particular, we benefit from local training processes (subdomain-dependent), which allows us to potentially avoid using the stochastic gradient method or$/$and improve its convergence. Typically, the mini-batches would actually correspond to training points for the local loss functions under consideration. At Schwarz iteration $k$, we hence minimize

\begin{eqnarray*}
\left.
\begin{array}{lcl}
\mathcal{L}^{\pm,(k)}({\boldsymbol w}) & = & \big\| \partial_t N^{\pm, (k)}({\boldsymbol w}^{\pm},\cdot,\cdot) -\nu({\boldsymbol x})\triangle N^{\pm, (k)}({\boldsymbol w}^{\pm},\cdot,\cdot) -  {\boldsymbol a}({\boldsymbol x}) \cdot \nabla N^{\pm, (k)}({\boldsymbol w}^{\pm},\cdot,\cdot)\\
& &  - r({\boldsymbol x})N^{\pm, (k)}({\boldsymbol w}^{\pm},\cdot,\cdot) \big\|_{L^2(\Omega_{\varepsilon}^{\pm}\times (0,T))}  +  \lambda_{\textrm{Ext}}^{\pm} \big\|N^{\pm, (k)}({\boldsymbol w}^{\pm},\cdot,\cdot)\big\|_{L^2(\Lambda_{\varepsilon}^{\pm}\times (0,T))}\\
& &  + \lambda_{\textrm{Int}}^{\pm} \big\|\mathcal{T}_{\pm} N^{\pm, (k)}({\boldsymbol w}^{\pm},\cdot,\cdot) - \mathcal{T}_{\pm} N^{\mp, (k-1)}(\overline{{\boldsymbol w}}^{\pm},\cdot,\cdot)\big\|_{L^2(\Gamma_{\varepsilon}^{\pm}\times (0,T))}  \, ,
\end{array}
\right.
\end{eqnarray*}
where $\overline{{\boldsymbol w}}^{\pm}$ were computed at the Schwarz iteration $k-1$. Recall that practically the loss functions are numerically evaluated by approximating the norm using training points, typically randomly chosen in $\Omega^{\pm}_{\varepsilon}$. This method allows for a complete spatial decoupling of the problem over 2 (or arbitrary number of) subdomains.  Finally, the reconstructed solution is hence defined as $N_{|\Omega_{\varepsilon}^{\pm}}(\cdot,t)=N^{\pm}(\cdot,t)$ for all $t \geq 0$. More specifically

\begin{eqnarray}\label{CVTOT2}
\lim_{k \rightarrow +\infty}\big\| \hspace{0.2cm}\|N^{+,(k)}(\overline{{\boldsymbol w}}^+,\cdot,\cdot)-N^{-,(k)}(\overline{{\boldsymbol w}}^+,\cdot,\cdot)\|_{\infty,\overline{\Omega}_{\varepsilon}^+\cap\overline{\Omega}_{\varepsilon}^-}\big\|_{L^{2}(0,T)}=0 \, ,
\end{eqnarray}
and we define the solution to the advection-diffusion-reaction equation as

\begin{eqnarray*}
N & = & 
\left\{
\begin{array}{ll}
N^{+,(k^{\textrm{cvg}})}(\overline{{\boldsymbol w}}^+,\cdot,\cdot), & \textrm{ in } \Omega_{\varepsilon}^+ \times (0,T) \, ,\\
N^{-,(k^{\textrm{cvg}})}(\overline{{\boldsymbol w}}^-,\cdot,\cdot), & \textrm{ in } \Omega_{\varepsilon}^- \times (0,T) \, .
\end{array}
\right.
\end{eqnarray*}

Practically, in order to evaluate the loss functions, it is necessary to compute the equation at some very large $N_TN_{\boldsymbol x}^{\pm}$ randomly chosen training points $\{({\boldsymbol x}^{\pm}_{j},t)\}_{n;j}$ in $\Omega_{\varepsilon}^{\pm}\times (0,T)$, as the $L^2-$norms are not exactly performed. From the point of view of the computation of the loss function (requiring the evaluation of the PDE at the training points), the algorithm is hence trivially embarrassingly parallel.  From the optimization point of view, the method now requires minimizing two loss functions. Naturally, the computation of the minimization problems is embarrassingly parallel as the two IBVPs are {\it totally} decoupled. As we are now considering two IBVPs on smaller spatial domains, we can locally adapt the depth of the local networks.\\
It is important to mention that, unlike SWR methods combined with standard numerical (finite-difference, -volume, -elements, pseudospectral) methods, for which convergence can be proven, the combination of SWR and PINN methods will not necessarily ensure convergence to zero of the residual history. This is due to the fact that from one Schwarz iteration to the next, the reconstructed solutions may slightly differ as the minima obtained by minimization of the local loss functions will a priori slightly differ. This fact is actually inherent to the NN-based method. However, we expect the residual history to be small from a practical point of view and for loss functions sufficiently small. In addition to this argument, let us mention that the transmission condition is naturally not exactly satisfied if it is included in the loss function. A large enough weight can, for instance, be imposed on the transmission constraint to ensure that it is accurately satisfied.


\subsection{About the interest of using SWR DDM for NN-based algorithms}\label{subsec:interest}
The estimation of the loss function using the direct PINN-method for solving local PDE is trivially embarrassingly parallel, as the estimation is independently performed at any given training point. However, this associated minimization problem (batch-size related) is not {\it locally specific}, and Stochastic Gradient Method (SGM) is hence an essential tool. In the proposed approach, the local loss functions which are evaluated have specific meanings; and allows to get accurate approximations of the solution in any given subdomain.

The SWR method is a domain decomposition method {\it in space} for solving PDE. Using standard advection-diffusion-reaction equation solvers, the main algorithmic costs are the loss function estimations and the computation of solutions to linear systems at each time iteration, involved in implicit or semi-implicit stable schemes \cite{halpern3}.  The latter has a polynomial complexity $O(N^{\alpha})$, where $1<\alpha\leq 3$ is typically dependent on the structure of the matrix involved in the linear system. Using a PINN approach, there are naturally no more linear systems to solve ``to estimate'' the solution. Instead, an optimization algorithm is necessary to parameterize the NN-based solution. Denoting by $N_{{\boldsymbol x};t}$ the (a priori very large) number of space-time training points to construct the loss function, and $N_W$ the total number of parameters. The computation of solution using the direct PINN method is decomposed into two parts:
\begin{itemize}
\item Estimation of the loss function, with a complexity $O(N_{{\boldsymbol x};t})$ with $N_{{\boldsymbol x};t}\gg 1$. This step is hence embarrassingly parallel for local PDE (with or without combination with the SWR method.)
\item Minimization of the loss function with a complexity $O(N_W^p)$ for $p \in (1,3)$. Typically stochastic gradient methods \cite{sgm0,sgm1,sgm2} are used to deal with possibly high dimensionality (for very accurate solutions) of this minimization problem and allows for a relatively  efficient  parallelization. 
\end{itemize}

Within the framework of DDM and for two subdomains, the SWR-NN indeed requires the embarrassingly parallel minimization of two independent loss functions constructed using local training. The empirical argument which justifies the proposed methodology is as follows. The structure and complexity of the solution is thought to be ``divided'' in the two (much more in practice of course) spatial subdomains. As a consequence, based on the local structure of the solutions in $\Omega^{\pm}_{\varepsilon}$, the depth of the local neural networks $N^{\pm}$ can then be adapted$/$reduced compared to the one-domain PINN approach with one unique neural network. The extreme case in that matter, would be a domain decomposition into small finite volumes, where the solution would be approximated by a constant (cell-center finite volume method) that is 0-depth NNs, even if the overall solution has a complex spatial structure. Naturally, the larger the subdomain size, the deeper the depth of the searched local neural network associated to this subdomain. For two subdomains, the minimization step within the SWR-NN consists in solving  {\it in parallel}
\begin{eqnarray*}
\overline{{\boldsymbol w}}^{\pm,(k)} & = & \textrm{argmin}_{{\boldsymbol w} \in W^{\pm}}\mathcal{L}^{\pm,(k)}({\boldsymbol w}) \, ,
\end{eqnarray*}
rather than (for direct method)
\begin{eqnarray*}
\overline{{\boldsymbol w}} & = & \textrm{argmin}_{{\boldsymbol w} \in W}\mathcal{L}({\boldsymbol w}) \, ,
\end{eqnarray*}
where $N_{W^{\pm}}:=\# W^{\pm} \leq N_{W}:=\# W$. That is, it is possible to decompose the minimization problem in several spatial subregions, where the spectral structure of the solution can be very different from on subdomain to the next, requiring locally smaller depths than using a unique global deep NN. Hence, for SWR-NN method we have to perform the following tasks.

\begin{itemize}
\item Estimation of the loss functions with a complexity $O(N^{\pm}_{{\boldsymbol x};T})$. This step is hence embarrassingly parallel within the SWR method and allows to deal with local (subdomain-dependent) training points.
\item Minimization of the local loss functions with a complexity $O(N_{W_+}^{\alpha})+ O(N^{\alpha}_{W_{-}})$ for $\alpha \in (1,3)$, where in principle $N_{W_{\pm}}\leq N_W$.
\end{itemize}
The downsize of SWR-NN methods is that it requires to repeat $k^{\textrm{cvg}}$ times (that until convergence) the computation of the uncoupled systems.
Unlike standard SWR-DDM, where the gain is on the computation of local linear systems of smaller size, the main interest is that we locally solve local (and less complex) minimization problems, where we expect the size of the {\it search space to be smaller}.

Notice that the SWR-DDM allows for an embarrassing parallelization of the overall PINN PDE solver. Indeed, unlike the standard computation of the (local) minima of the loss function, which requires non-trivial non-embarrassingly parallelization, the proposed approach allows for the embarrassingly parallel computation of minima of local loss functions. Three levels of parallelization are then possible
\begin{itemize}
\item Trivial parallelization of the estimation of the local loss functions.
\item Embarrassingly parallel computation of the minima of the local loss functions.
\item In addition, the minimum of a local loss function can also be performed in parallel using the domain decomposition method for PINN, as proposed in see \cite{pinns4}.
\end{itemize}
From the computational point of view, the SWR-PINN algorithm allows i) to adapt the depth of (most of) local NNs compared to using one unique (global) NN, and ii) to estimate the local loss functions using local subdomain-dependent training points and potentially allows for using direct  (none-stochastic) gradient methods for a sufficiently large number of subdomains. This step is the analog of the reduction of the size of the linear systems to be solved (scaling effect) within standard SWR when are applied as real space solvers \cite{XA1,XA2,XA3}.\\
\\
We here summarize the overall computational complexity of SWR-PINN and direct PINN methods.
\begin{itemize}
\item Direct approach: $O(N_{{\boldsymbol x;t}})+O(N^P_W)$. In this case, we expect $N_W$ to be large, and $P$ depends on the used optimization algorithm.
\item SWR approach: $O\big(k^{\textrm{cvg}}(N^P_{W_+}N^{P}_{W_-})\big)+O(N^-_{{\boldsymbol x;t}})+O(N^+_{{\boldsymbol x;t}})$. In this case, we expect $N_{W^{\pm}} < N_{W}$. As $P$ is strictly greater than $1$ there is a scaling effect which makes this approach potential more efficient. Moreover the prefactor is also though to be much smaller using SWR methods. Practically it is required for $k^{\textrm{cvg}}$ to be small enough. As it is well-known the choice of the transmission conditions is a crucial element to minimize $k^{\textrm{cvg}}$. Dirichlet transmission conditions is known to provide very slow convergence. At the opposite of the spectrum and for wave-like equations, Dirichlet-to-Neumann like transmission conditions are known to provide extremely fast convergence, but can be computationally complex to approximate. Another way to accelerate the convergence of the SWR algorithm consists in increasing the subdomain overlap (that is increase $\varepsilon$) For the advection-diffusion-reaction equation, optimized SWR method, based on optimized Robin-transmission conditions is a good compromise between convergence rate and computational complexity \cite{halpern2}. As specified above, the computation of the loss function is embarrassingly parallel unlike the direct approach. 
\end{itemize}

\subsection{Nonlocal operator}\label{subsec:nonlocal}
We have argued above that the use of SWR methods allows for an efficient parallel computation of the overall loss functions through the efficient estimation (using local training points) of local loss functions. We show below that whenever nonlocal terms are present in the equation, the efficiency of the SWR-PINN method is not deteriorated by those terms.  In the following, we assume that the equation contains a nonlocal operator $F$, typically defined as a convolution product:
\begin{itemize}
\item  $F(u) = \partial_x^{\alpha}u$, with $\alpha>0$, fractional derivative in space, modeling nonlocal effect \cite{whatis}. The latter is actually defined as a convolution. We refer to \cite{review} for details.
\item  $F(u) = \rho*_{\boldsymbol x}u$ where $*_{\boldsymbol x}$ denotes the spatial convolution product, a nonlocal potential for some given function $\rho$.
\end{itemize}
We consider the equation on a truncated domain $\Omega \subset \R^d$ with boundary $\Gamma$, as follows
\begin{eqnarray}
\label{e1bis}
\left\{
\begin{array}{l}
\partial_t u = \nu({\boldsymbol x})\triangle u + F(u), \, {\boldsymbol x} \in  \Omega, \, t>0,\\
u({\boldsymbol x},0) = u_0({\boldsymbol x}), \, {\boldsymbol x} \in  \Omega ,\\
u = 0, \, {\boldsymbol x} \in  \Gamma \,  ,
\end{array}
\right.
\end{eqnarray}
and such that $F$ is defined as a convolution product in space
\begin{eqnarray*}
F(u)  = u*_{\boldsymbol x}\rho  =  \int_{\Omega} u({\boldsymbol x}-{\boldsymbol y},t)\rho({\boldsymbol y})d{\boldsymbol y} \, ,
\end{eqnarray*}
Then the SWR-PINN scheme reads
\begin{eqnarray}\label{owrs3}
\left\{
\begin{array}{lcl}
\partial_t N^{\pm, (k)} & = & \nu({\boldsymbol x})\triangle N^{\pm, (k)} +  \int_{\Omega_{\varepsilon}^{\pm}}N^{\pm,(k)}({\boldsymbol w}^{\pm},{\boldsymbol x}-{\boldsymbol y},t)\big)\rho({\boldsymbol y})d{\boldsymbol y} \\
& & + \int_{\Omega_{\varepsilon}^{\mp}}N^{\mp,(k-1)}(\overline{{\boldsymbol w}}^{\mp},{\boldsymbol x}-{\boldsymbol y},t)\big)\rho({\boldsymbol y})d{\boldsymbol y} , \, \hbox{ in } \Omega_{\varepsilon}^{\pm} \times (0,T),\\
N^{\pm, (k)}(\cdot,0) & =& u_0^{\pm},\, \hbox{ in } \Omega_{\varepsilon}^{\pm}, \\
\mathcal{T}_{\pm} N^{\pm, (k)} & = & \mathcal{T}_{\pm} N^{\mp, (k-1)}, \, \hbox{ on } \Gamma_{\varepsilon}^{\pm} \times (0,T) ,\\
N^{\pm, (k)} & = & 0, \, \hbox{ on } \Lambda \backslash\Gamma_{\varepsilon}^{\pm} \times (0,T) \, ,
\end{array}
\right.
\end{eqnarray}
where $\overline{\boldsymbol {w}}^{\pm}\in W_{\pm}$ was computed at the previous Schwarz iteration, with some transmission operator $\mathcal{T}_{\pm}$.  Hence in this case, we still have to minimize local loss functions

\begin{eqnarray*}
\left.
\begin{array}{lcl}
\mathcal{L}^{\pm,(k)}({\boldsymbol w}) & = & \big\|\partial_t N^{\pm, (k)}({\boldsymbol w}^{\pm},\cdot,\cdot) -\nu({\boldsymbol x})\triangle N^{\pm, (k)}({\boldsymbol w}^{\pm},\cdot,\cdot)\\
& &  -\int_{\Omega_{\varepsilon}^{\pm}}N^{\pm,(k)}({\boldsymbol w}^{\pm},{\boldsymbol x}-{\boldsymbol y},t)\big)\rho({\boldsymbol y})d{\boldsymbol y} \\
& & - \int_{\Omega_{\varepsilon}^{\mp}}N^{\mp,(k-1)}(\overline{{\boldsymbol w}}^{\mp},{\boldsymbol x}-{\boldsymbol y},t)\big)\rho({\boldsymbol y})d{\boldsymbol y}\big\|_{L^2(\Omega_{\varepsilon}^{\pm}\times (0,T))} \\
& &  + \lambda_{\textrm{Int}}^{\pm} \big\|\mathcal{T}_{\pm} N^{\pm, (k)}({\boldsymbol w}^{\pm},\cdot,\cdot) - \mathcal{T}_{\pm} N^{\mp, (k-1)}(\overline{{\boldsymbol w}}^{\mp},\cdot,\cdot)\big\|_{L^2(\Gamma_{\varepsilon}^{\pm}\times (0,T))} \\
& & +  \lambda_{\textrm{Ext}}^{\pm} \big\|N^{\pm, (k)}({\boldsymbol w}^{\pm},\cdot,\cdot)\big\|_{L^2(\Lambda_{\varepsilon}^{\pm}\times (0,T))} \, .
\end{array}
\right.
\end{eqnarray*}

Practically, we can approximate the convolution product as follows. Denoting by $\{{\boldsymbol x}_j\}_{j\in \mathcal{J}^{\pm}}\in \Omega^{\pm}$ the local spatial training points, for ${\boldsymbol x}\in \Omega_{\varepsilon}^{\pm}$, we approximate
\begin{eqnarray*}
\int_{\Omega_{\varepsilon}^{\pm}}N^{\pm,(k)}({\boldsymbol w}^{\pm},{\boldsymbol x}-{\boldsymbol y},t)\big)\rho({\boldsymbol y})d{\boldsymbol y} +  \int_{\Omega_{\varepsilon}^{\mp}}N^{\mp,(k-1)}(\overline{{\boldsymbol w}}^{\mp},{\boldsymbol x}-{\boldsymbol y},t)\big)\rho({\boldsymbol y})d{\boldsymbol y}
\end{eqnarray*}
by
\begin{eqnarray*}
\left.
\begin{array}{l}
 \sum_{j\in \mathcal{J}^{\pm}}c_jF\big(N^{\pm,(k)}({\boldsymbol w}^{\pm},{\boldsymbol x}_i-{\boldsymbol x}_j,t)\big) + \sum_{j\in \mathcal{J}^{\mp}}c_jF\big(N^{\mp,(k-1)}(\overline{{\boldsymbol w}}^{\mp},{\boldsymbol x}_i-{\boldsymbol x}_j,t)\big)  \, , 
\end{array}
\right.
\end{eqnarray*}

for some weights $\{c_j\}_{j\in \mathcal{J}^{\pm}}$. \\
\\
 As it was discussed above the interest of using a DDM is to decompose the training and search of the local solution over smaller set of parameters. However, whenever the equation is {\it nonlocal}, it is necessary to extend the search of the parameters in the global computational domains. More specifically, for local equations, the local NN-solution in $\Omega_{\varepsilon}^+$ (resp. $\Omega_{\varepsilon}^-$) only requires parameters in $W_+$ (resp. $W_-$). However, if the equation is nonlocal, in order to construct the solution in $\Omega_{\varepsilon}^{\pm}$, we have in principle to search the NN parameters in all $W$, containing both $W_+$ and $W_-$, as the solution values in $\Omega_{\varepsilon}^{\pm}$ depend on values of the solution in $\Omega^{\mp}_{\varepsilon}$. This problem would also occur to construct the loss function, using the direct PINNs method, within the term
\begin{eqnarray*}
\int_{\Omega}N({\boldsymbol w},{\boldsymbol x}-{\boldsymbol y},t)\big)\rho({\boldsymbol y})d{\boldsymbol y} \, .
\end{eqnarray*}

The SWR-PINN method allows to deal with this issue, as at Schwarz iteration $k$, the loss function in $\Omega^{\pm}_{\varepsilon}$ is evaluated through the solution in $\Omega_{\varepsilon}^{\pm}$ at the previous Schwarz iteration $(k-1)$ from $N^{\mp}(\overline{{\boldsymbol w}}^{\mp},{\boldsymbol x},t)$ thanks to the previously evaluated parameter $\overline{{\boldsymbol w}}^{\mp}$.

\subsection{How about a non-iterative domain decomposition in space?}\label{subsec:direct}
The domain decomposition method which is derived in this paper is a SWR-in space method which is an iterative method allowing for the convergence of the decomposed solution towards the exact solution of the PDE under consideration. The main weakness of this DDM is the fact that the decoupled system has to be solved several times (iterative method). It is hence natural to ask if a ``similar spatial domain decomposition'', but non-iterative, is possible.\\
In this goal, we decompose the domain $\Omega$ as above: $\Omega_{\varepsilon}^{\pm}$, with or without overlap ($\Omega_{\varepsilon}^+\cap\Omega_{\varepsilon}^-=\emptyset$ or $\Omega_{\varepsilon}^+\cap\Omega_{\varepsilon}^-\neq \emptyset$), with $\epsilon \geq 0$ and consider \eqref{e1bis}. That is we search for a solution of the form $N_{|\Omega_{\varepsilon}^{\pm}}(\cdot,t)=N^{\pm}(\overline{{\boldsymbol w}}^{\pm}\cdot,t)$ for all $t\geq 0$ such that
\begin{eqnarray}\label{owrs4}
\left\{
\begin{array}{lcl}
\partial_t N^{\pm} & = & \nu({\boldsymbol x})\triangle N^{\pm} +  \int_{\Omega_{\varepsilon}^{\pm}}N^{\pm}({\boldsymbol w}^{\pm},{\boldsymbol x}-{\boldsymbol y},t)\big)\rho({\boldsymbol y})d{\boldsymbol y} \\
& & + \int_{\Omega_{\varepsilon}^{\mp}}N^{\mp}({\boldsymbol w}^{\mp},{\boldsymbol x}-{\boldsymbol y},t)\big)\rho({\boldsymbol y})d{\boldsymbol y} , \, \hbox{ in } \Omega_{\varepsilon}^{\pm} \times (0,T),\\
N^{\pm} & = & 0, \, \hbox{ on } \Gamma \times (0,T) \, ,
\end{array}
\right.
\end{eqnarray}
where $\overline{\boldsymbol {w}}^{\pm}\in W_{\pm}$.  The PINN-method then consists in solving

\begin{eqnarray*}
\overline{{\boldsymbol w}}^{\pm} & = & \textrm{argmin}_{{\boldsymbol w} \in W^{\pm}}\mathcal{L}^{\pm}({\boldsymbol w}) \, ,
\end{eqnarray*}
where 
\begin{eqnarray*}
\left.
\begin{array}{lcl}
\mathcal{L}^{\pm}({\boldsymbol w}) & = & \big\| \partial_t N^{\pm, (k)}({\boldsymbol w}^{\pm},\cdot,\cdot) -\nu({\boldsymbol x})\triangle N^{\pm, (k)}({\boldsymbol w}^{\pm},\cdot,\cdot)\\
& &  -\int_{\Omega_{\varepsilon}^{\pm}}N^{\pm}({\boldsymbol w}^{\pm},{\boldsymbol x}-{\boldsymbol y},t)\big)\rho({\boldsymbol y})d{\boldsymbol y} \\
& & - \int_{\Omega_{\varepsilon}^{\mp}}N^{\mp}({\boldsymbol w}^{\mp},{\boldsymbol x}-{\boldsymbol y},t)\big)\rho({\boldsymbol y})d{\boldsymbol y}\big\|_{L^2(\Omega_{\varepsilon}^{\pm}\times (0,T))} \\
& & +  \lambda_{\textrm{Ext}}^{\pm} \big\|N^{\pm, (k)}({\boldsymbol w}^{\pm},\cdot,\cdot)\big\|_{L^2(\Gamma\times (0,T))} \, .
\end{array}
\right.
\end{eqnarray*}

Therefore, in this case, we still have to minimize local loss functions. However, there are 2 main issues:
\begin{itemize}
\item Even if $\varepsilon$ is taken equal to zero, the decoupling of the solution in the two subdomains, naturally induces a discontinuity at the subdomain interfaces. It is possible to impose additional {\it compatibility conditions} (to be included in the loss function), in the form of continuity condition $N^{+}=N^{-}$ at $\Gamma_{\varepsilon}^{\pm}$, differentiability, but the reconstructed global solution $N$ (such that $N_{|\Omega^{\pm}_{\varepsilon}}=N^{\pm}$), will obviously not be an approximate solution to the equation under consideration. Moreover, the compatibility conditions will induce a re-coupling of the two systems in the spirit of the following item.

\item The two systems, in $N^+$ and $N^-$ are actually also coupled through the nonlocal term. This effect is similar to the addition of a compatibility condition described above. Hence, the computation of the loss functions $\mathcal{L}^{\pm}$ would not be embarrassingly parallel anymore. This is not an issue in the SWR framework; as in the latter case, say at Schwarz iteration $k$, the nonlocal term uses the approximate solution at the Schwarz iteration $k-1$, which is a known quantity.
\end{itemize} 
Hence unlike the SWR-PINN method, for which \eqref{CVTOT2} occurs 
\begin{eqnarray}\label{CVTOT3}
N^{+}_{\overline{\Omega}_{\varepsilon}^+\cap\overline{\Omega}_{\varepsilon}^-}(\overline{{\boldsymbol w}}^+,\cdot,\cdot) \neq N_{\overline{\Omega}_{\varepsilon}^+\cap\overline{\Omega}_{\varepsilon}^-}^{-}(\overline{{\boldsymbol w}}^-,\cdot,\cdot) \, .
\end{eqnarray}

\section{Numerical Experiments}\label{sec:numerics}
\begin{figure}[hbt!]
\begin{center}
\includegraphics[height=6cm,keepaspectratio]{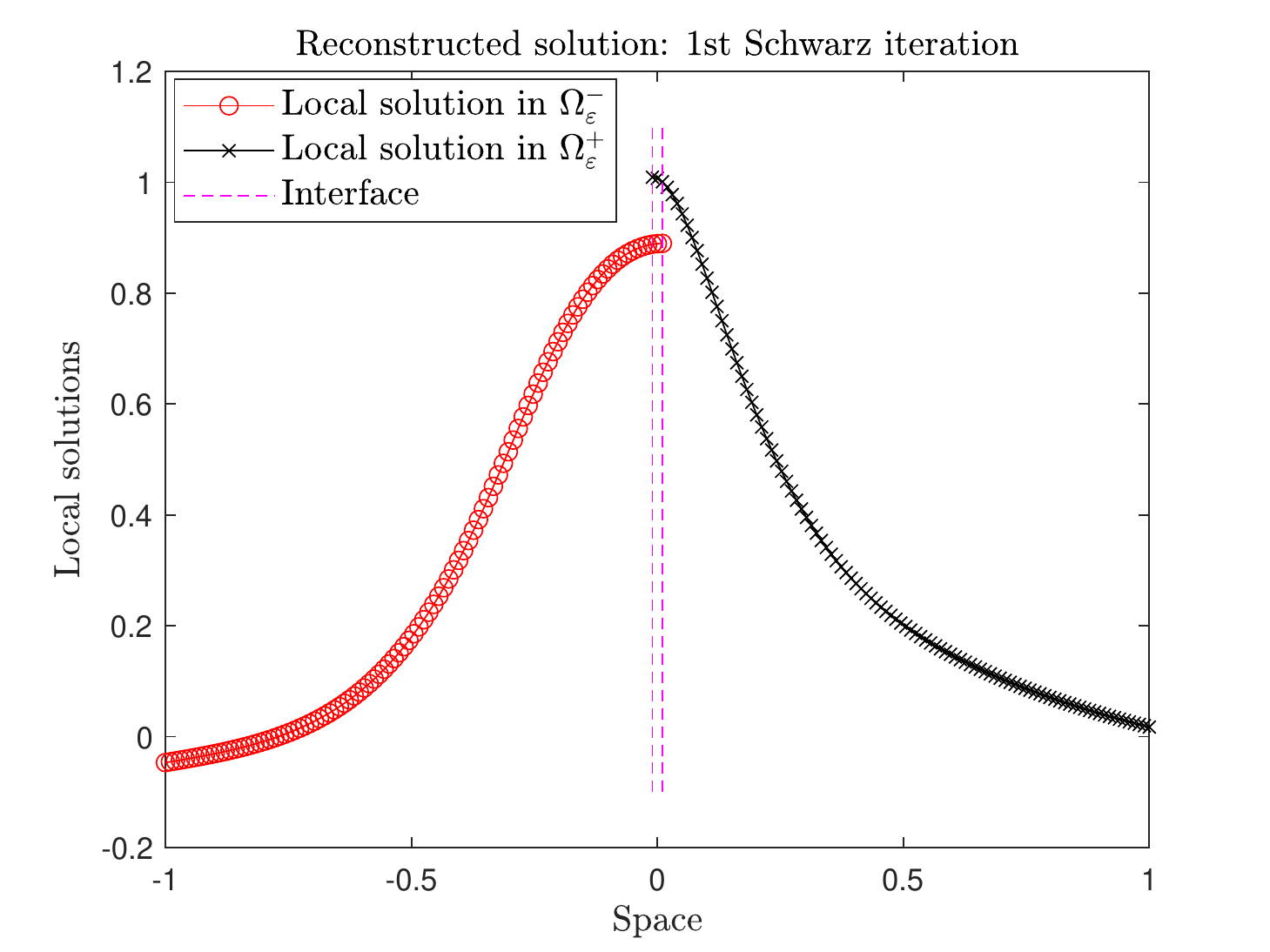}
\includegraphics[height=6cm,keepaspectratio]{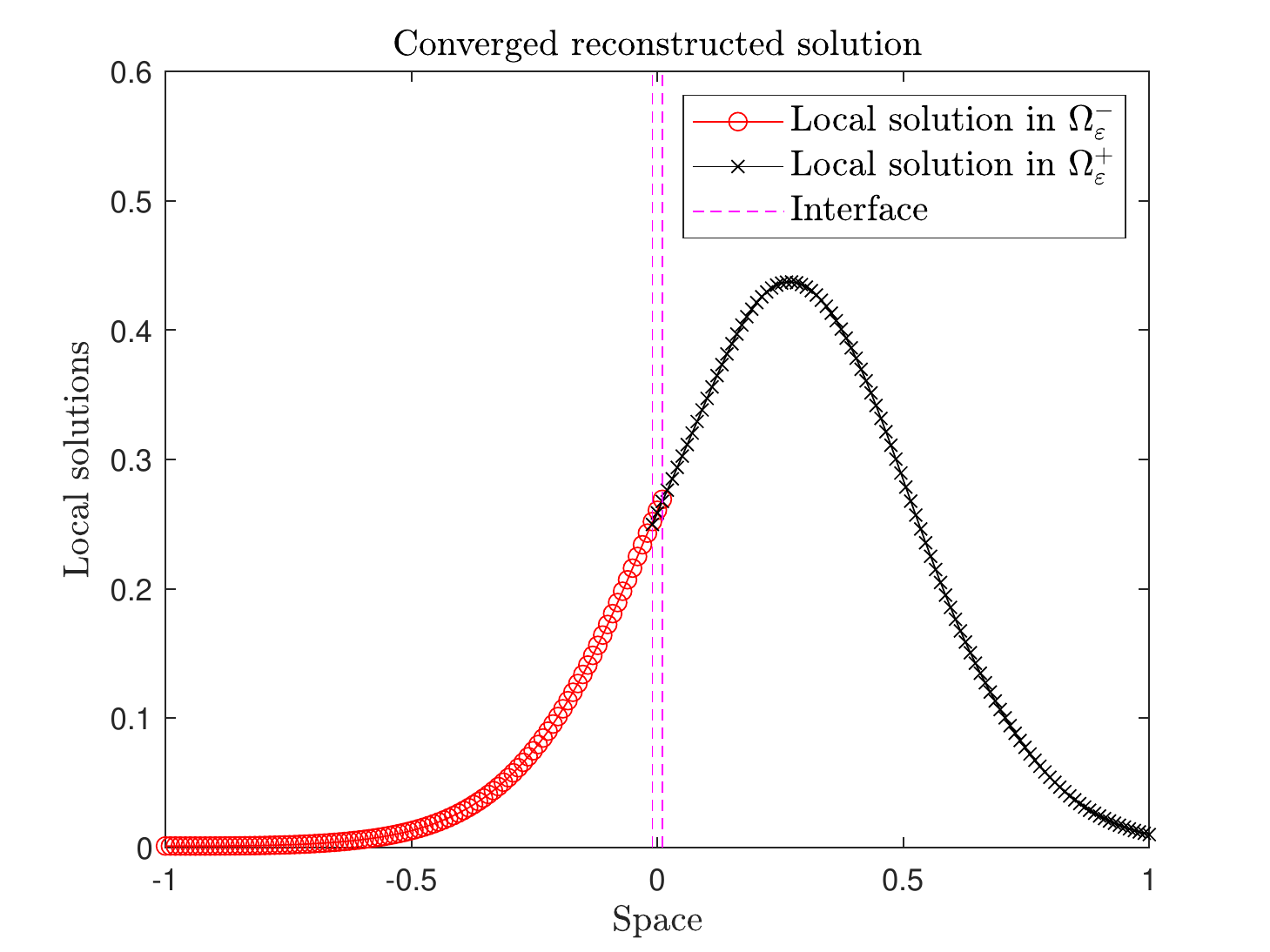}
\end{center}
\caption{{\bf Experiment 1.} (Left) Reconstructed solutions after the first Schwartz iteration (Right) Reconstructed solutions at convergence of the SWR method.}
\label{fig1}
\end{figure}

In this section, we propose basic experiments in order to numerically illustrate the convergence of the overall method.  The PINN-algorithm was implemented using deep-learning and optimization toolboxes from {\tt matlab}, {\tt DeepXDE} \cite{deepXDE} and {\tt tensorflow} \cite{tensorflow2015-whitepaper}. Although relatively simple, these experiments illustrate the proof of concept of the proposed strategy and not to provide the best convergence possible (which will be the purpose of a future work). \\
\noindent{\bf Experiment 1.} We consider the standard advection-diffusion-reaction equation 
\begin{eqnarray}\label{adr}
\partial_t u & = &  a \partial_x u + \nu \partial_{xx}u + r u \, ,
\end{eqnarray}
on $\Omega\times[0,T]=(-1,1)\times [0,0.25]$ with Dirichlet boundary conditions at $x=\pm 1$ and such that $a=0.5$, $\nu=0.0$,  $r=0.1$ and the initial conditions is $u_0(x)=\exp(-30(x-0.1)^2)$. We decompose the domain in two subdomains $\Omega^-_{\varepsilon}=(-1,\varepsilon/2)$ and $\Omega^+_{\varepsilon}=(-\varepsilon/2,1)$ with $\varepsilon=0.1$. We here use the Classical Schwarz Waveform Relaxation method, based on Dirichlet transmission conditions $N^{\pm,(k)}(\cdot,\pm\varepsilon/2,\cdot) = N^{\mp,(k-1)}(\cdot,\mp \varepsilon/2,\cdot)$, where $N^{\pm}$ are the two local NN defined in $\Omega_{\varepsilon}^{\pm}$. We consider the following data: the NN have both $5$ layers, with $25$ neurons each. We select $5000$ internal collocation points. We also use a local SGM with $5$ epochs and mini-batches size of $500$. In the gradient method the learning rate with decay rate of $5\times 10^{-3}$ starting at $10^{-3}$.  We reconstruct the overall solution using a total of $5000$ prediction points. Initially we take $N^{\pm,(0)}=0$.  We report the reconstructed solution after the first SWR iteration (resp. converged SWR algorithm) in Fig. \ref{fig1} (Left) (resp. \ref{fig1} (Right)) from two the local solutions in $\Omega_{\varepsilon}^{\pm}$ at final time $T=0.25$, and overlapping zone of size $1/99$.  

\begin{figure}[hbt!]
\begin{center}
\includegraphics[height=6cm,keepaspectratio]{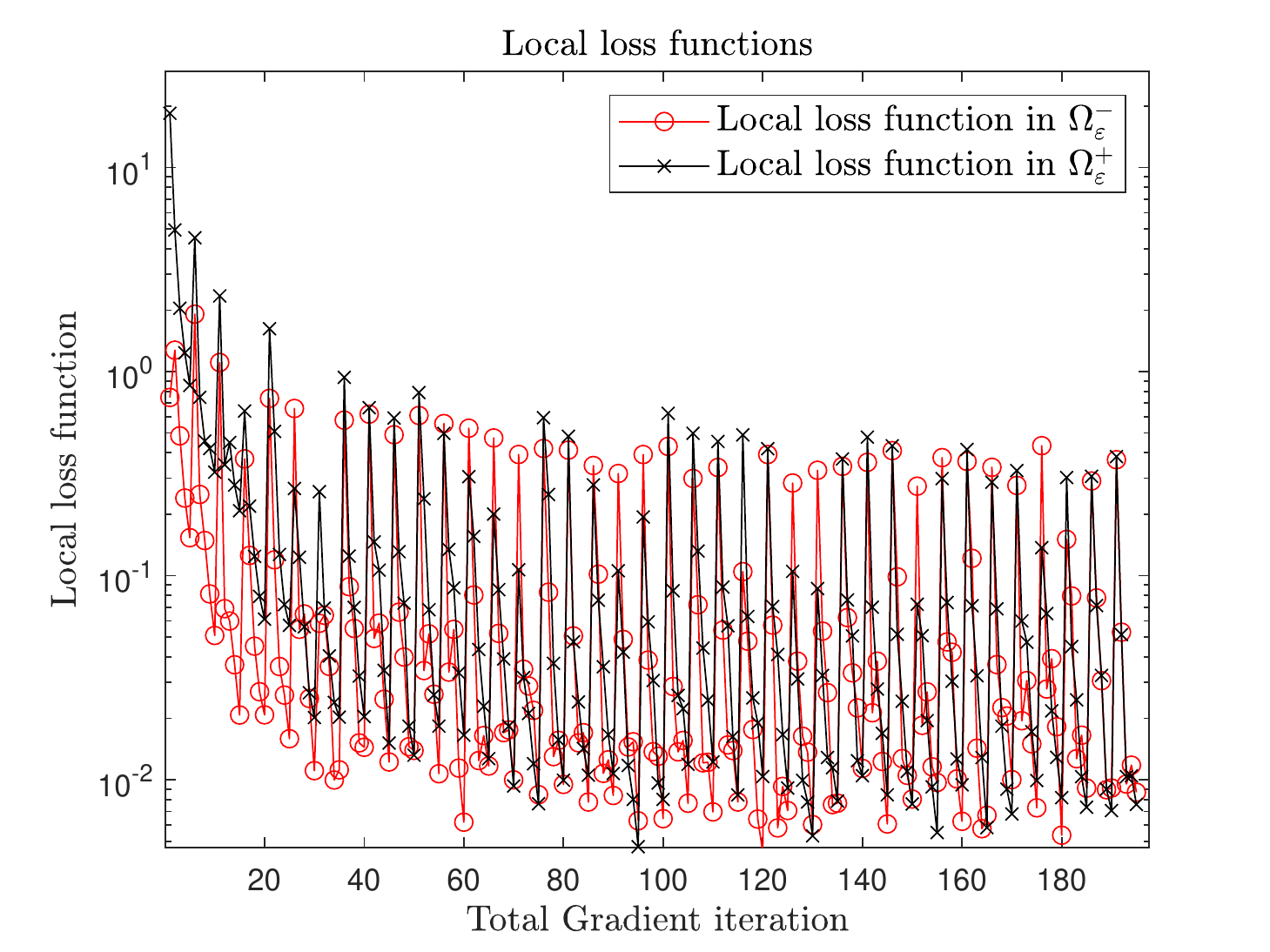}
\includegraphics[height=6cm,keepaspectratio]{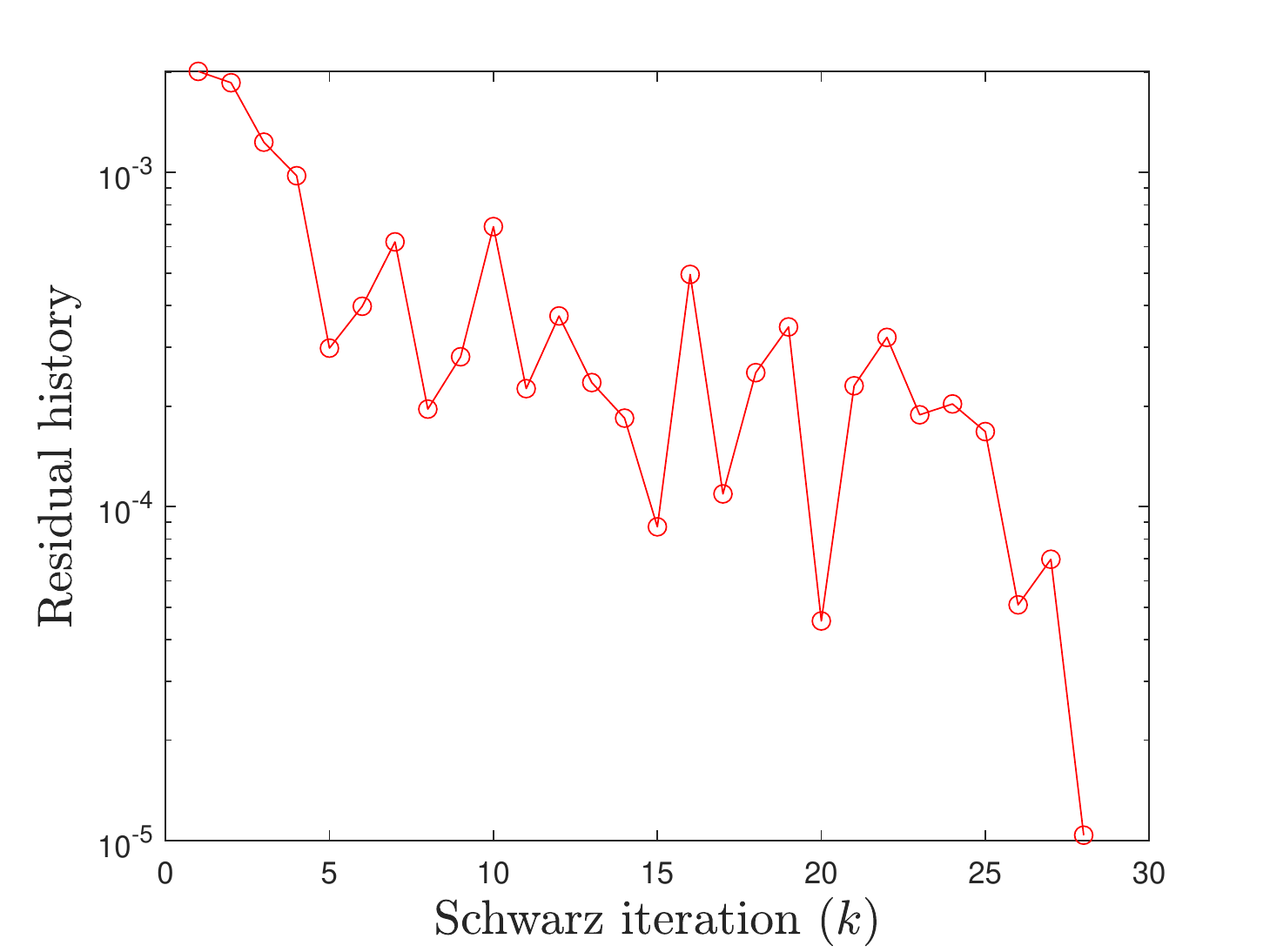}
\end{center}
\caption{{\bf Experiment 1.} (Left) Converged reconstructed solutions (Right) Graph of local loss functions.}
\label{fig2}
\end{figure}

 The SWR convergence rate is defined  as the slope of the logarithm of the residual history according to  the Schwarz iteration number, that is $\{(k,\log(\mathcal{E}^{(k)})) \, : \, k \in \N\}$, with (for 2 subdomains)
\begin{eqnarray}\label{rate2}
\mathcal{E}^{(k)} := \sum_{i=1}^{2}\big\| \hspace{0.2cm}\|N^{+,(k)}-N^{-,(k)}\|_{\infty;\overline{\Omega}_{\varepsilon}^+\cap\overline{\Omega}_{\varepsilon}^-}\big\|_{L^{2}(0,T)}\leq  \delta^{\textrm{Sc}},
\end{eqnarray}
 $\delta^{\textrm{Sc}}$ being a small parameter.\\
We report in Fig \ref{fig2} (Left) the graph of convergence of the stochastic gradient methods applied to each local loss functions. Notice that each ``oscillation'' corresponds to a new Schwarz iteration. We report in Fig. \ref{fig2} (Right) the graph of convergence of the SWR-method in the form of the residual history in the overlapping zone. The combination of convergent SWR methods with standard numerical (finite element, finite-difference) methods for which there is a uniform convergence to zero of the residual history as a function of the Schwarz iterations. 

\begin{figure}[hbt!]
\begin{center}
\includegraphics[height=6cm,keepaspectratio]{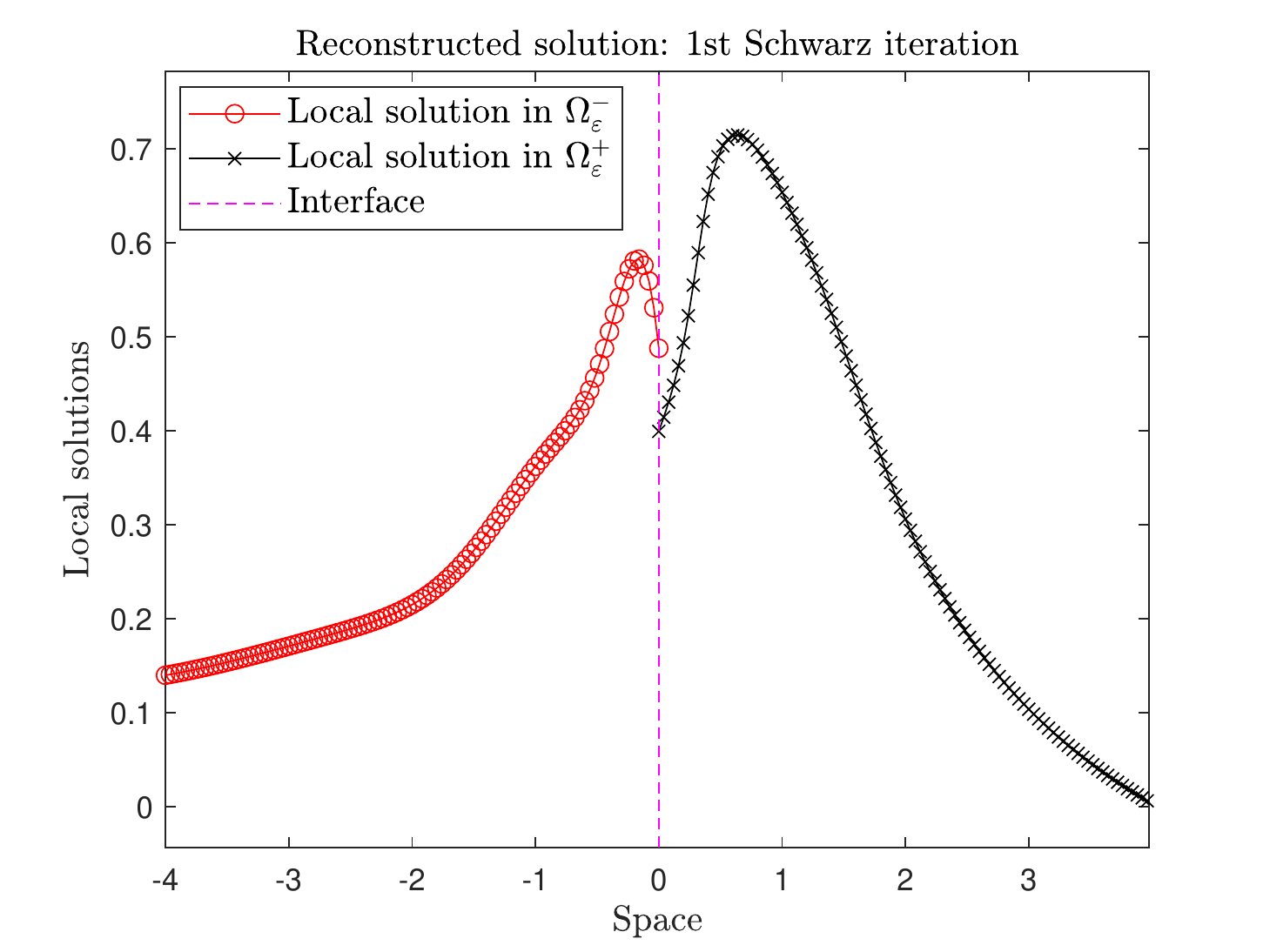}
\includegraphics[height=6cm,keepaspectratio]{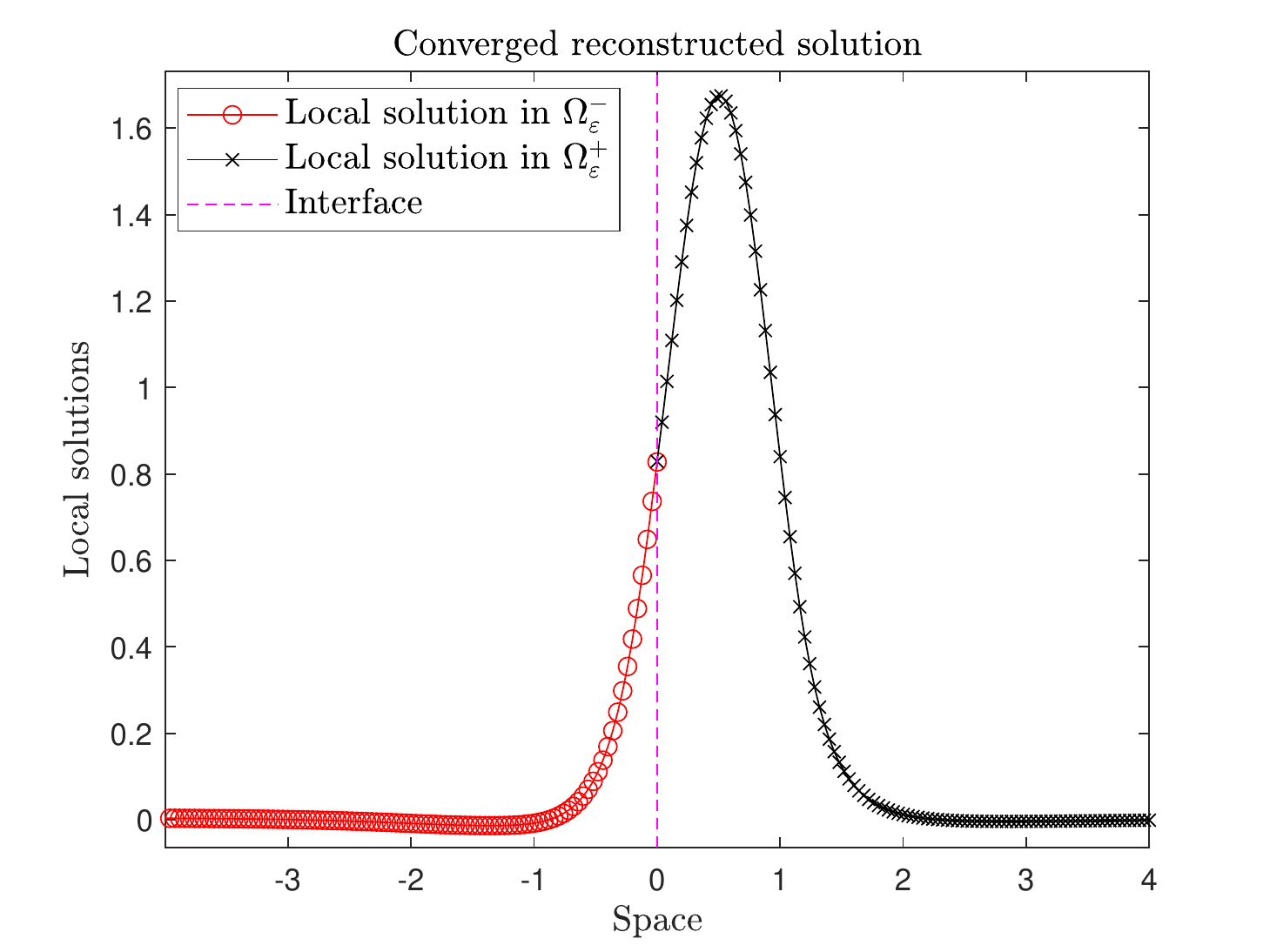}
\end{center}
\caption{{\bf Experiment 2.} (Left) Reconstructed solutions after the first Schwartz iteration (Right) Reconstructed solutions at convergence of the SWR method.}
\label{fig3}
\end{figure}

Notice that, unlike the converged SWR method combined with standard numerical methods where the residual goes to zero when $k$ goes to infinity, the residual does not exactly go to zero. This is due to the fact that from one Schwarz iteration to the next, the (local) solution are obtained by constructing ``new'' local minima as the local loss functions are small but {\it not} null, and hence change from one iteration to the next.  
\\
\\
{\bf Experiment 2.} In the following experiment, we implement a Robin-SWR method for solving \eqref{adr}, which is expected to provide better convergence than CSWR \cite{halpern2}. As it was discussed in \cite{halpern2} and recalled above, the optimized SWR (and more generally Robin-SWR) methods is convergent, even without overlap, that is when $\varepsilon$ is null. We consider the same equation as above with $a=1$, $r=0$, $\nu=5\times 10^{-2}$ and $T=0.5$. The initial conditions is $u_0(x)=10\exp(-4x^2) + $. We decompose the domain in two subdomains $\Omega_{\varepsilon}^-=(-4,0)$ and $\Omega_{\varepsilon}^+=(0,4)$, with hence $\varepsilon=0$. The Robin transmission conditions $(\lambda \pm \partial_x)N^{\pm,(k)}(\cdot,0,\cdot) = (\lambda + \partial_x) N^{\mp,(k-1)}(\cdot,0,\cdot)$, where $N^{\pm}$ are the two local NN defined in $\Omega_{\varepsilon}^{\pm}$ and we have taken $\lambda=5$. We consider the following data: the NN have both $8$ layers, with $30$ neurons each. We select $5000$ internal collocation points. We also use local SGM with $5$ epochs and mini-batches size of $500$. In the gradient method the learning rate with decay rate of $10^{-3}$ starting at $10^{-3}$.  We reconstruct the overall solution using a total of $5000$ prediction points.  Initially we take $N^{\pm,(0)}=0$.  We report the reconstructed solution after the first SWR iteration (resp. converged SWR algorithm) in Fig. \ref{fig3} (Left) (resp. \ref{fig3} (Right)) from two the local solutions in $\Omega_{\varepsilon}^{\pm}$ at final time $T=0.5$.

\begin{figure}[hbt!]
\begin{center}
\includegraphics[height=6cm,keepaspectratio]{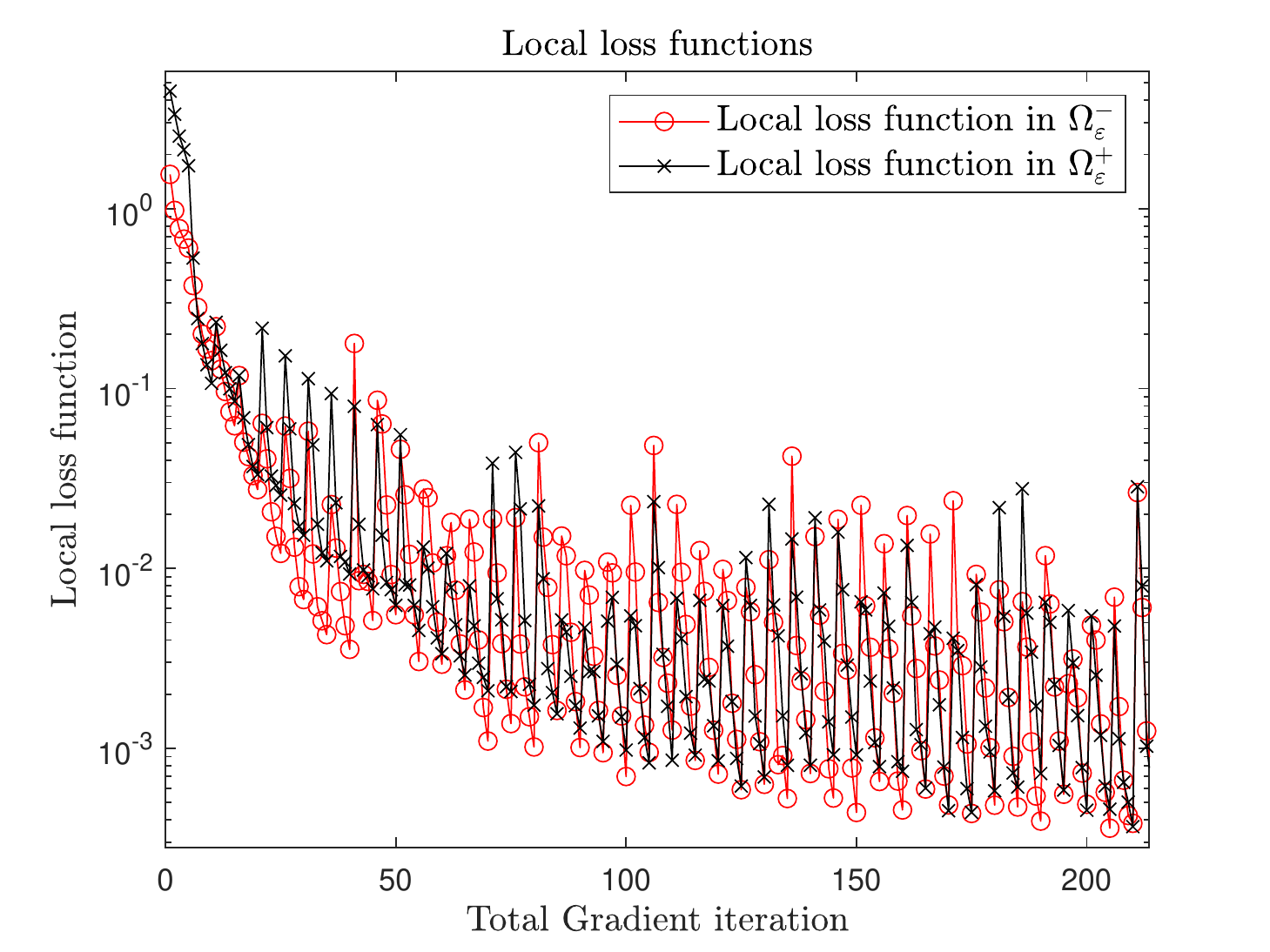}
\includegraphics[height=6cm,keepaspectratio]{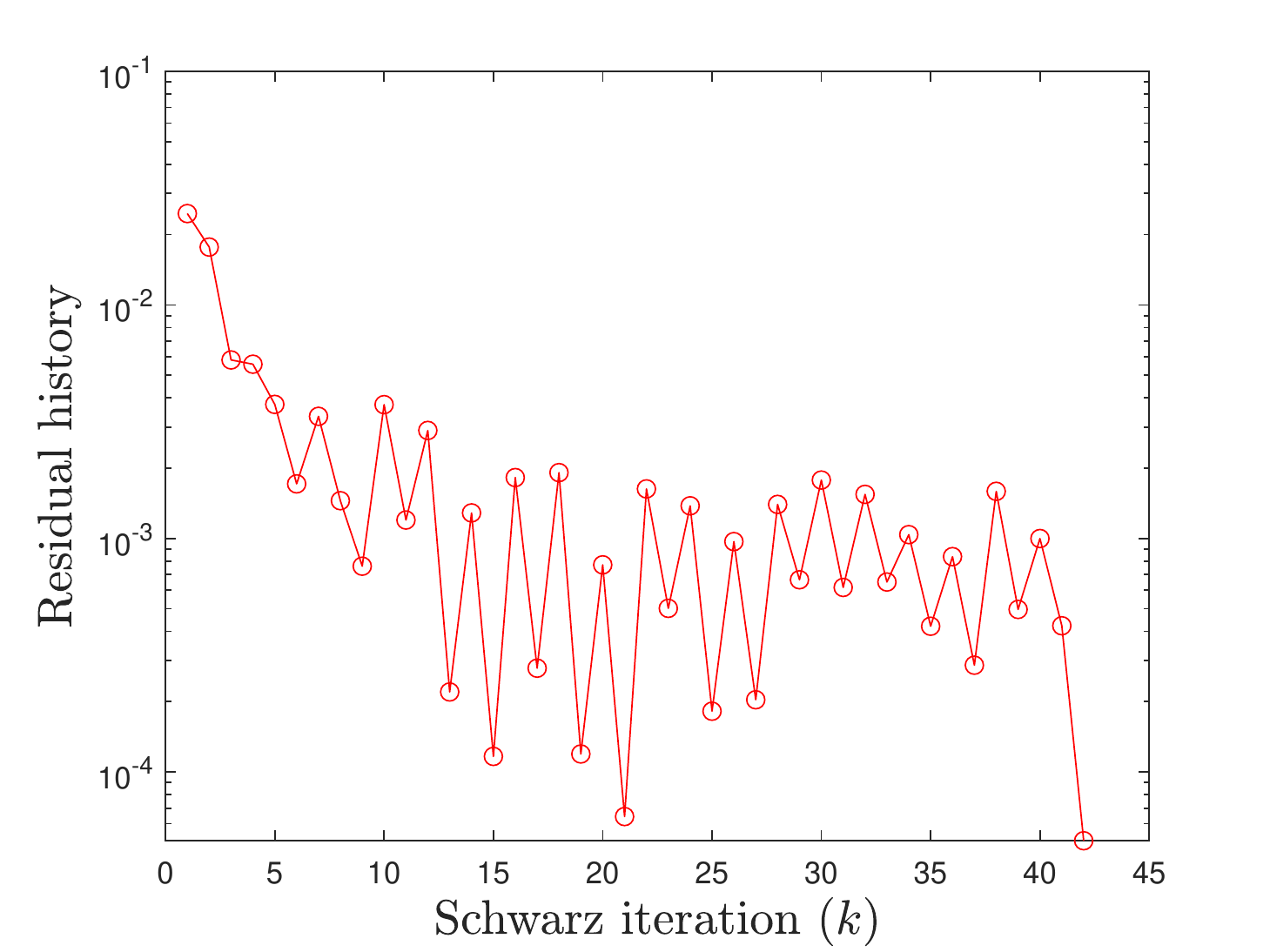}
\end{center}
\caption{{\bf Experiment 2.} (Left) Converged reconstructed solutions (Right) Graph of local loss functions.}
\label{fig4}
\end{figure}

We next report in Fig \ref{fig4} (Left) the graph of convergence of the stochastic gradient methods applied to each local loss functions.  We report in Fig. \ref{fig4} (Right) the graph of convergence of the SWR-method in the form of the residual history in the overlapping zone.
\\
\\
 Importantly, we observe that Robin-SWR-PINN still converges even if the two subdomains do not overlap.

\begin{figure}[hbt!]
\begin{center}
\includegraphics[height=6cm,keepaspectratio]{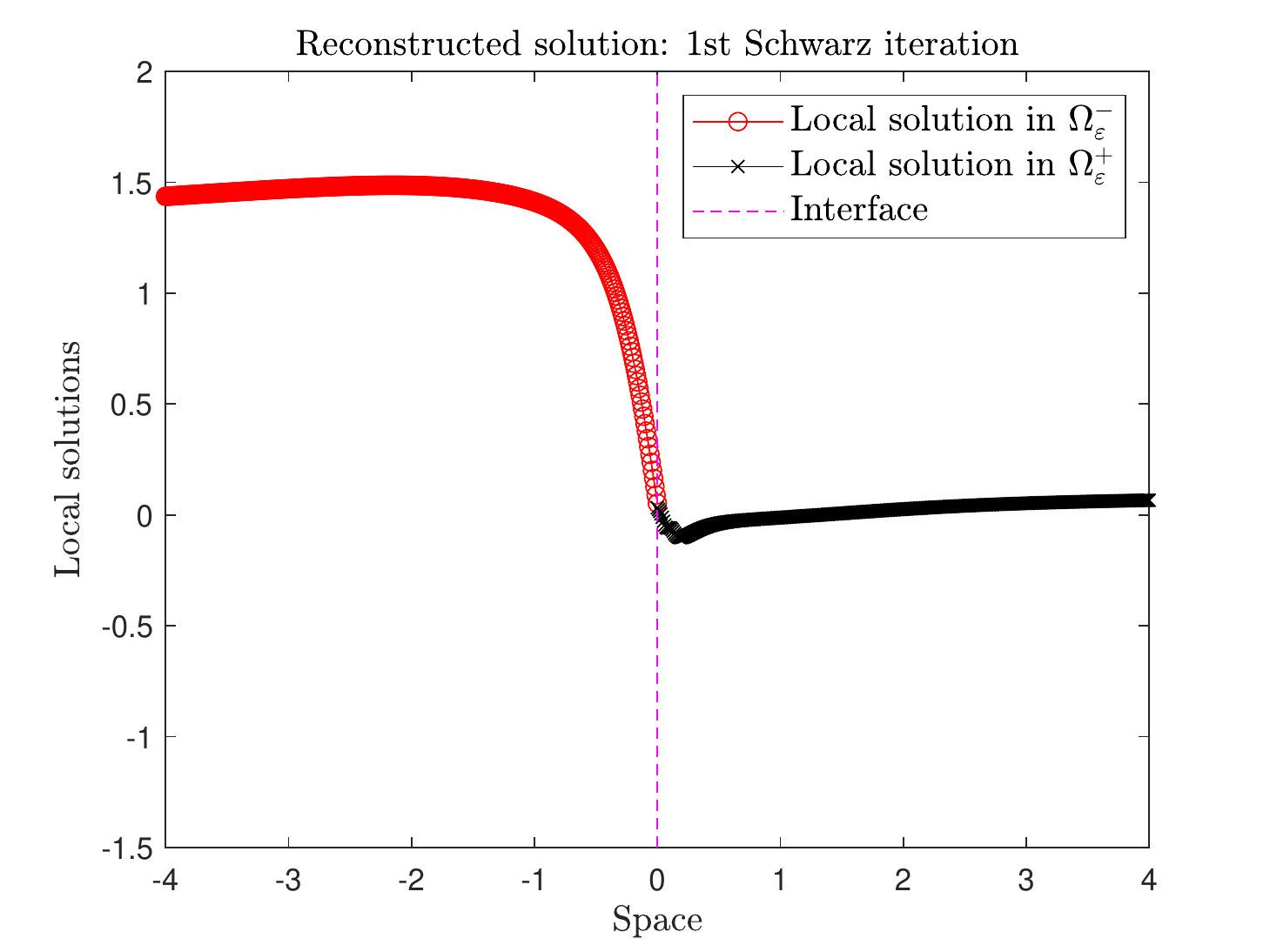}
\includegraphics[height=6cm,keepaspectratio]{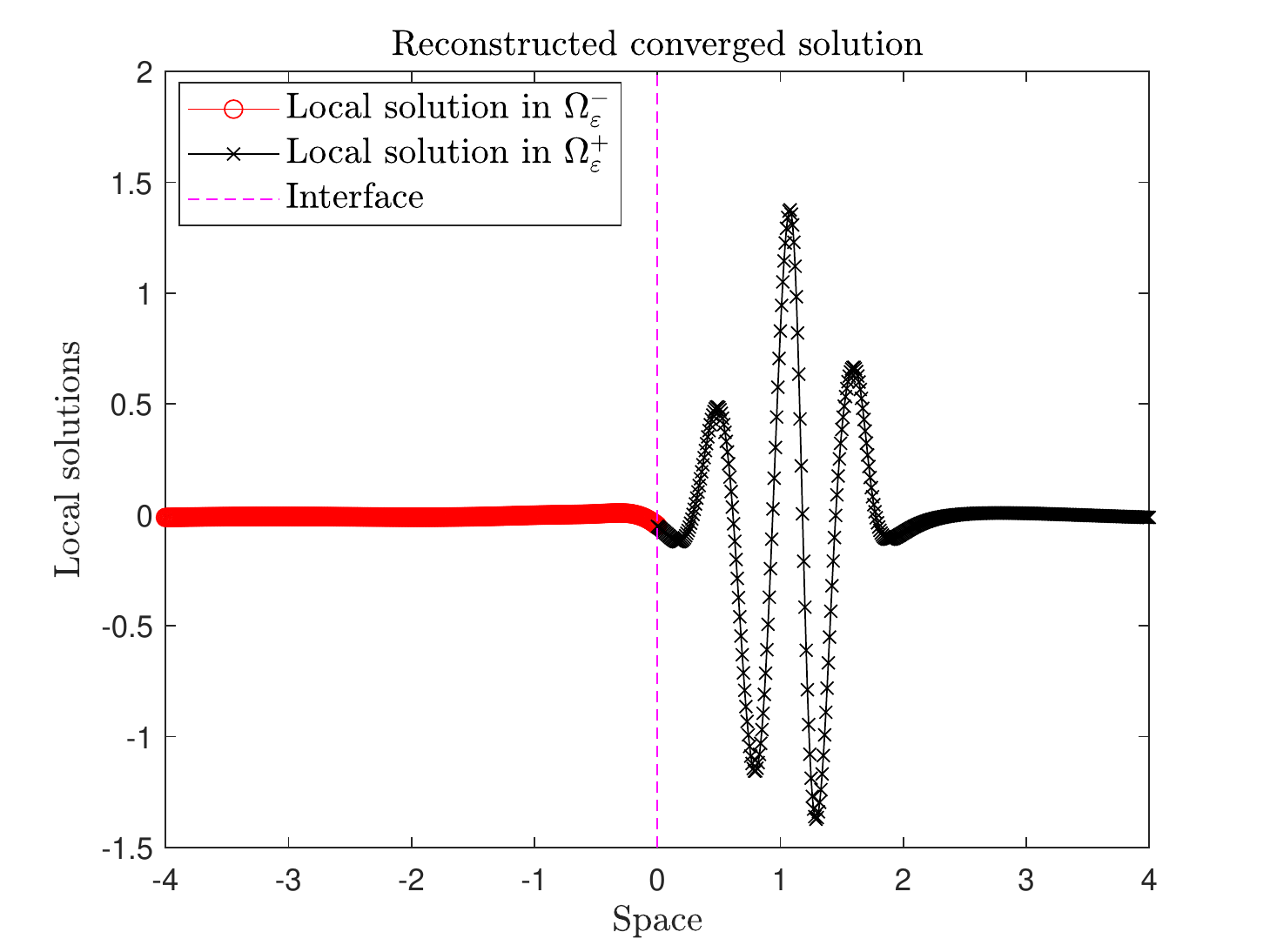}
\includegraphics[height=6cm,keepaspectratio]{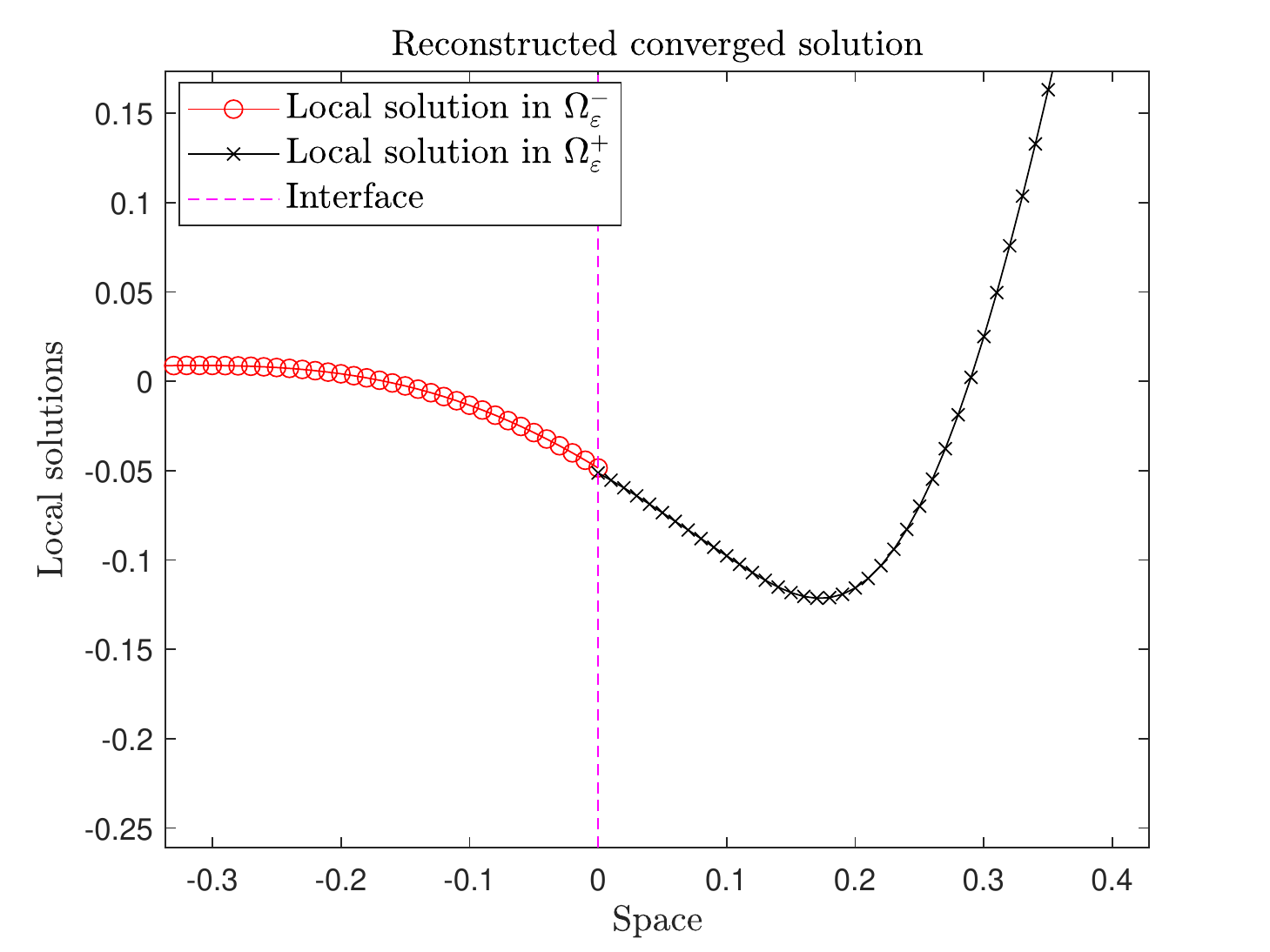}
\includegraphics[height=6cm,keepaspectratio]{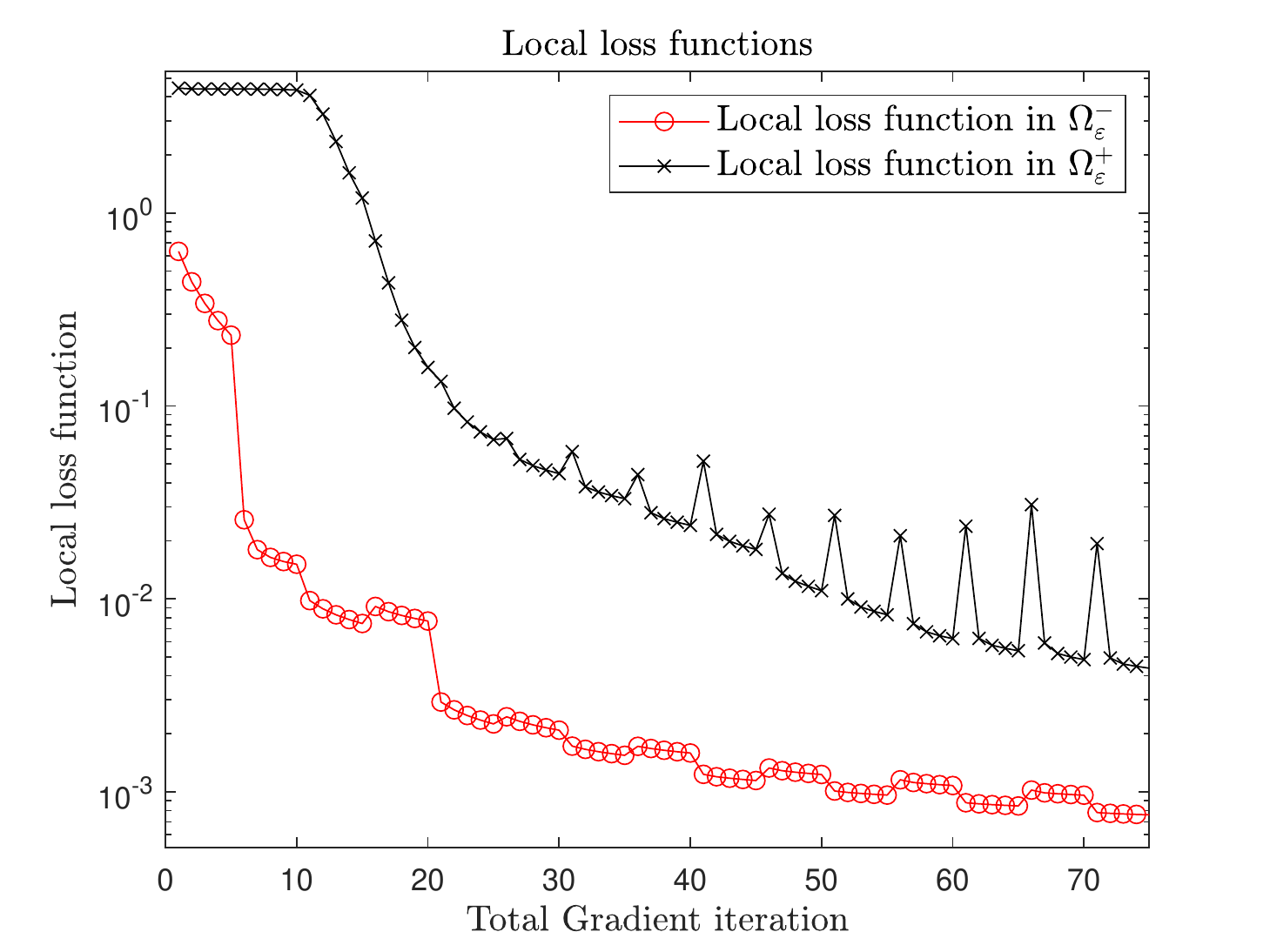}
\end{center}
\caption{{\bf Experiment 2bis.} (Top-Left) Converged solutions after the first Schwarz iteration. (Top-Right) Reconstructed converged solution. (Bottom-Left) Zoom in of the converged solution in the interface region. (Bottom-Right) Local loss function convergence.}
\label{fig5}
\end{figure}

{\bf Experiment 2bis.}
In the following non-overlapping 2-domain Robin-SWR experiment, we now consider that the diffusion coefficient is {\it space-dependent}; more specifically, $\nu(x)=0.1$ (resp. $2.5\times 10^{-2}$) for $x \in \Omega_{\varepsilon}^{+}$ (resp. $x \in \Omega_{\varepsilon}^{-}$). The rest of the data are as follows: $a=1$, $r=0$, $\nu=5\times 10^{-2}$, $T=0.1$, $\lambda=5$. The initial condition is given by $u_0(x)=10\exp(-3(x-1)^2)\cos(10(x-1)^2)$ and is such that the solution has a very different structure in the two subdomains. We want here to illustrate the ability of the derived approach to select different depths of the local neural networks, depending on the structure of the solution: in $\Omega^-_{\varepsilon}$ (resp. $\Omega_{\varepsilon}^+$) the solution is mainly null (resp. oscillatory), except close to the interface. The two subdomains are $\Omega_{\varepsilon}^-=(-2,0)$ and $\Omega_{\varepsilon}^+=(0,2)$, with hence $\varepsilon=0$. The two local NN $N^{\pm}$, over $\Omega_{\varepsilon}^{\pm}$ have the following structure: $N^{-}$ (resp. $N^+$) possesses $3$ (resp. $10$) layers and $10$ (resp. $50$) neurons. {\it The minimization process in $N^-$ is much more efficiently performed than in $N^+$ with a relatively similar accuracy}. As above, we select $5000$ internal collocation points. We also use local SGM with $5$ epochs and mini-batches size of $5000$. In the gradient method the learning rate with decay rate of $10^{-3}$ starting at $10^{-3}$.  We reconstruct the overall solution using a total of $20000$ prediction points.  Initially we take $N^{\pm,(0)}=0$.  We report the reconstructed solution after the first SWR iteration (resp. converged SWR algorithm) in Fig. \ref{fig5} (Top-Left) (resp. \ref{fig5} (Top-Right)) from two the local solutions in $\Omega_{\varepsilon}^{\pm}$ at final time $T=0.1$. We also zoom in (\ref{fig5}, (Bottom-Left)), in the interface region to better observe the SWR convergence. The local loss functions are represented in Fig. \ref{fig5} (Bottom-Right). We observe that roughly, the computational time to perform the solution in $\Omega_{\varepsilon}^{+}$ was $2.5$ times faster than in $\Omega_{\varepsilon}^{-}$.

\begin{figure}[hbt!]
\begin{center}
\includegraphics[height=6cm,keepaspectratio]{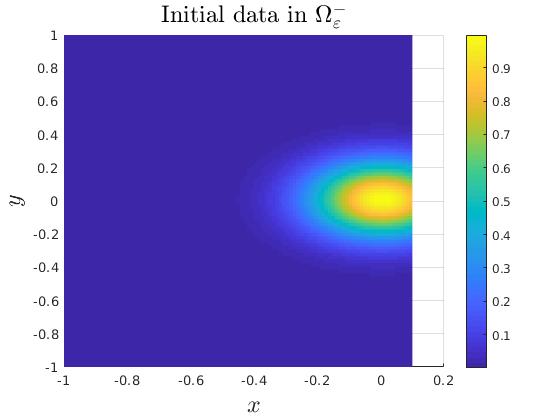}
\includegraphics[height=6cm,keepaspectratio]{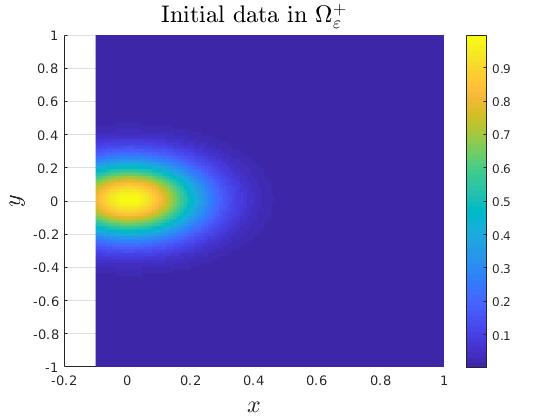}
\includegraphics[height=5.85cm,keepaspectratio]{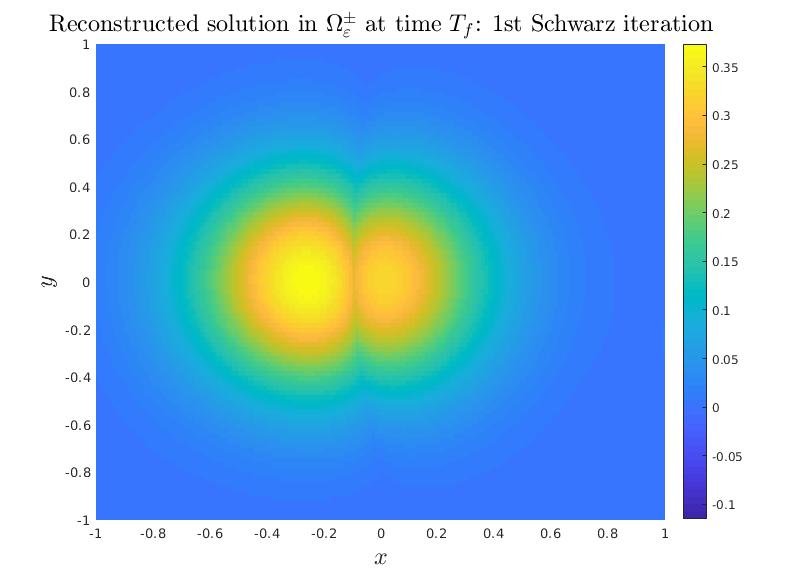}
\includegraphics[height=5.85cm,keepaspectratio]{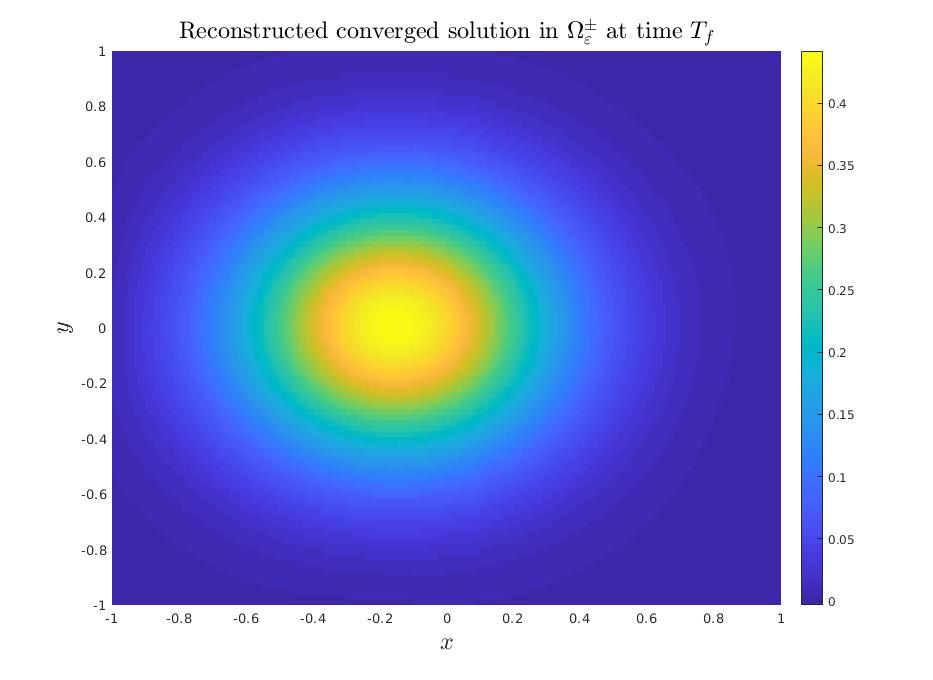}
\end{center}
\caption{{\bf Experiment 3.} (Top) Initial data in $\Omega_{\varepsilon}^{\pm}$. (Bottom-Left) Solutions after the first Schwarz iteration at time $T_f$ . (Bottom-Right) Reconstructed converged solution at time $T_f$.}
\label{fig6}
\end{figure}

\noindent{\bf Experiment 3.}  In this last experiment, we consider a two-dimensional advection-diffusion equation on a square $[-1,1]^2\times [0,T]$.
\begin{eqnarray*}
\partial_t u +{\bf a}\cdot \nabla u -\nu \triangle u & = & 0 \, ,
\end{eqnarray*}
with ${\bf a}=(-0.5,0)^T$ and $\nu=0.1$ and $T=2 \pi/10$.  The two subdomains are $\Omega_{\varepsilon}^{+}=(-1,\varepsilon)\times (-1,1)$ and  $\Omega_{\varepsilon}^{+}=(-1,\varepsilon)\times (-1,1)$ where $\varepsilon=0.1$; hence the interfaces are located at $\pm 0.1$. The initial data is a Gaussian function $u_0({\boldsymbol x})=\exp(-5\|{\boldsymbol x}\|^2)$ and the final computational time is $T=0.1$. A classical SWR algorithm is here combined with the PINN method.  On the other subdomain boundaries, we impose null Dirichlet boundary conditions. The equation is solved using the library {\tt DeepXDE} \cite{deepXDE} combined with {\tt tensorflow} \cite{tensorflow2015-whitepaper}. In each subdomain is used a neural network with $3$ layers and $10$ neurons; Adam's optimizer is used (learning rate $10^{-3}$, epoch=$10^3$) along with $\tanh$ activation function. In Fig. \ref{fig6} (Top), we report the initial data in $\Omega_{\varepsilon}^{\pm}$. In Fig. \ref{fig6}, we represent the solution at the end of the first Schwarz iteration (Left) and fully converged solution (Right) at final time $T_f$. In future works, we will propose more advanced simulations. The corresponding code is available on {\tt GitHub}, where the interested reader could find all the relevant information regarding the code.

\section{Conclusion}\label{sec:conclusion}

In this paper, we have derived a Schwarz Waveform Relaxed Physics-Informed Neural Networks (SWR-PINN) method for solving advection-diffusion-reaction equations in parallel. Some preliminary illustrating experiments are presented to validate the approach.

\subsection{Pros. and cons. of the SWR-PINN method}\label{subsec:pros-cons}
We summarize below the pros and cons of the proposed method.\\
\noindent{\bf Pros.}
\begin{itemize} 
\item Embarrassingly parallelization of the local loss function training.
\item Parallel construction of local neural networks with adaptive depth and complexity.
\item For convergent PINN algorithms, the SWR-PINN is convergent.
\item Flexible choice of the transmission conditions.
\end{itemize}
{\bf Cons.}
\begin{itemize}
\item As a fixed point method, SWR methods require several iterations.
\item The transmission conditions must be accurately satisfied through a penalization term in the loss function in order to accurately implement the SWR algorithm. Ideally, we should directly include the transmission within the definition in the NN.  This is possible, considering the CSWR (Dirichlet-based transmission conditions) method and the following NN, $T^{\pm}$
\begin{eqnarray*}
T^{\pm;(k)}({\boldsymbol w},{\boldsymbol x},t) & = & N_{|\Gamma_{\varepsilon}^{\pm}\times (0,T)}^{\mp;(k-1)}(\overline{{\boldsymbol w}}^{\mp},{\boldsymbol x},t) + \big({\bf 1}_{\Gamma_{\varepsilon}^{\pm}\times (0,T)}({\boldsymbol x},t)-1\big)N^{\pm;(k)}({\boldsymbol w},{\boldsymbol x},t) \, .
\end{eqnarray*}

\item Convergence or high precision of the overall algorithm can be hard to reach if the PINN algorithm is not used with sufficiently high precision. Instable numerical behavior can also be observed with the CSWR method.
\end{itemize}
\subsection{Concluding remarks and future investigations}
As far as we know, this paper is the first attempt to combine the SWR and PINN methods. Although the theory of SWR-DDM is now well developed in terms of convergence and convergence rate for different types of evolution PDE and their approximation with finite difference and finite element methods, the theory of convergence of PINN is not yet complete. Consequently, the convergence of the overall SWR-PINN method is still subject to the proof of convergence of the latter, which is largely empirically established. In future works, we plan to focus on ``real-life'' experiments where the main benefits of the SWR-PINN will be exhibited and illustrated.

\bibliographystyle{unsrt}
\bibliography{refs}

\end{document}